\documentclass[10pt,a4paper]{amsart}

\input{ort_coord.sty}
\input{ivan_style.sty}

\newcommand{\thi}{1.2mm}
\newcommand{\cthi}{.5mm}
\newcommand{\rad}{.2cm}

\begin{document}	
\setstcolor{red}
\title[Links in orthoplicial Apollonian packings]
{Links in orthoplicial Apollonian packings} 

\thanks{\noindent$^\dagger$ Partially supported by grant IEA-CNRS\\
\indent$^*$ Partially supported by the Austrian Science Fund (FWF), projects F-5503 and P-34763, and CNRS}
\author[Jorge L. Ram\'irez Alfons\'in]{Jorge L. Ram\'irez Alfons\'in$^\dagger$}
\address{
IMAG, Univ.\ Montpellier, CNRS, Montpellier, France }
\email{jorge.ramirez-alfonsin@umontpellier.fr}
\author[Iván Rasskin]{Iván Rasskin$^*$}
\address{LIS, Aix-Marseille Université, CNRS, Marseille, France}
\email{ivan.rasskin@lis-lab.fr}



\begin{abstract} 
In this paper, we establish a connection between Apollonian packings and knot theory. We introduce new representations of links realized in the tangency graph of the regular crystallographic sphere packings. Particularly, we prove that any algebraic link can be realized in the cubic section of the orthoplicial Apollonian packing. We use these representations to improve the upper bound on the ball number of an infinite family of alternating algebraic links. Furthermore, the later allow us to reinterpret the correspondence of rational tangles and rational numbers and to reveal geometrically primitive solutions for the Diophantine equation $x^4 + y^4 + z^4 = 2t^2$.
\end{abstract}

\maketitle
\section{Introduction}
Apollonian packings and their generalizations appear in many fields of science, including the modeling of granular systems \cite{Anishchik1995APOgranular}, fluid emulsions \cite{kwok2020APOemulsions} as well as in number theory \cite{apoNumber}. In this paper, we further explore the applications of Apollonian packings by introducing new representations of links realized in the tangency graph of certain three-dimensional analogues of Apollonian packings, extending thus their applications into the novel area of knot theory.

\subsection{Main results}
We begin by proving that any link can be realized in the tangency graph of any regular crystallographic sphere packing (Theorem \ref{lem:ortholinks}). We then focus our attention on algebraic links and show that any algebraic link can be realized in the tangency graph of a cubic Apollonian section of the orthoplicial Apollonian packing (Theorem \ref{thm:orthocubicrep}). The diagrams arising from this construction, called \textit{orthocubic representations}, have the following applications:

\subsubsection{Ball number}
A \textit{necklace representation} of a link $L$ is a sphere packing that contains a collection of disjoint cycles in its tangency graph realizing $L$. Necklace representations have been used to study the cusp volume of hyperbolic $3$-manifolds \cite{gabai2021hyperbolic}. The \textit{ball number} of $L$, denoted by $\mathrm{ball}(L)$, is defined as the minimum number of spheres needed to construct a necklace representation of $L$. It is known that $\mathrm{ball}(2_1^2)=8$ and $9\le\mathrm{ball}(3_1)\leq 12$ \cite{maehara2007}, where $2_1^2$ denotes the Hopf link and $3_1$ the trefoil knot. Currently, the Hopf link remains the only link for which the ball number is known. In \cite{RR20}, the authors provided a linear upper bound on the ball number of every nontrivial and nonsplittable link $L$, in terms of its \textit{crossing number}, denoted by $\mathrm{cr}(L)$, which is defined as the minimal number of crossings among all the diagrams of $L$. Specifically, it was shown that  $\mathrm{ball}(L)\le5\mathrm{cr}(L)$  and proposed the following.
\begin{conj}\label{conj:4n} For any nontrivial and nonsplittable link $L$, $\mathrm{ball}(L)\le 4\mathrm{cr}(L).$ Moreover, the equality holds if $L$ is alternating.
\end{conj}

\begin{figure}[H]
	\centering
	\includegraphics[width=.475\textwidth,align=t]{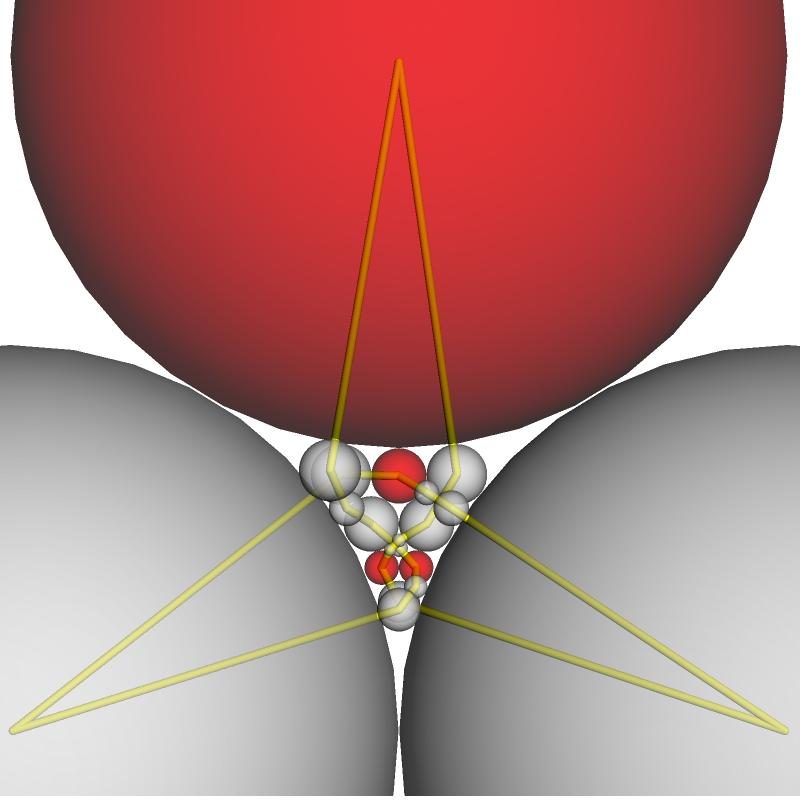}\hspace{.2cm}
	\includegraphics[width=.475\textwidth,align=t]{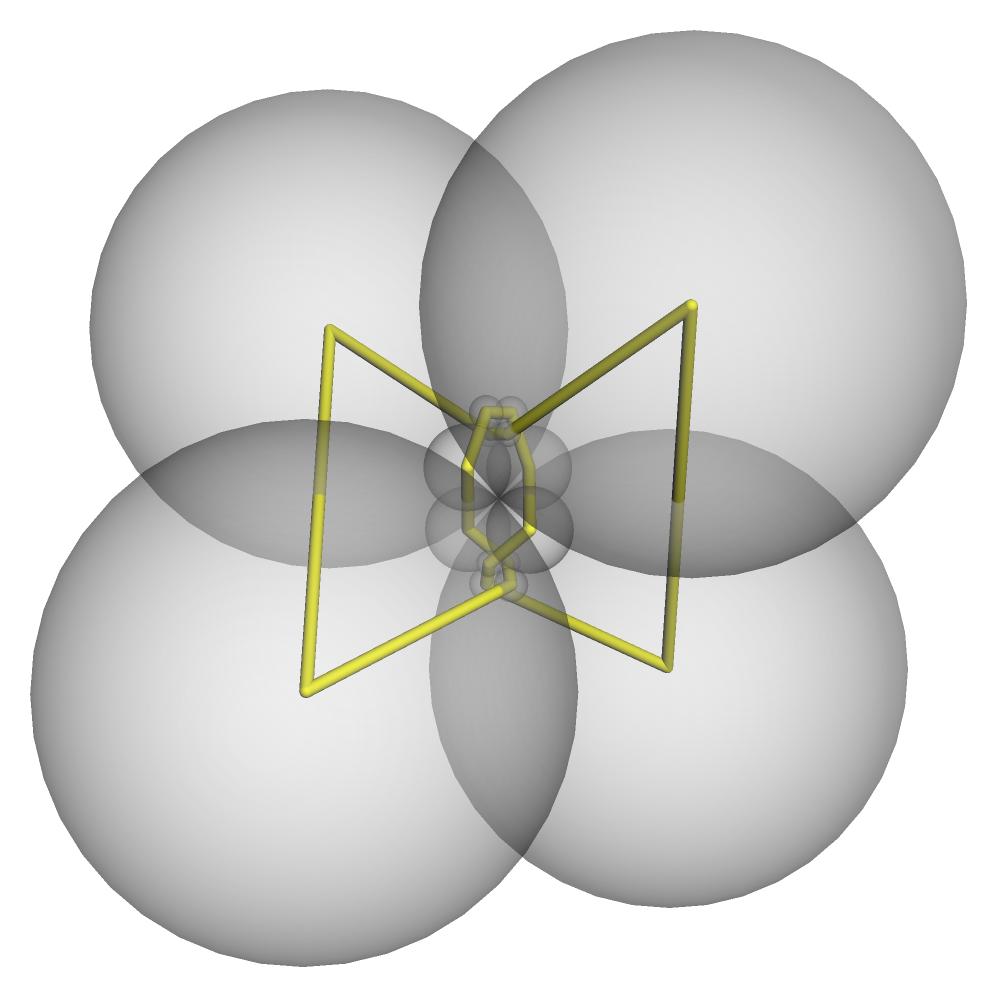}
	
	\caption{(Left) A necklace representation of the ``Figure-Eight" knot $4_1$  with 20 spheres obtained using the method from \cite{RR20}, featuring 4 connecting spheres (highlighted in red); (right) an orthocubic representation of the same knot with 16 spheres, without connecting spheres.}
	\label{fig:neck41}
\end{figure}

The inequality of Conjecture \ref{conj:4n} is motivated by the observation that, locally, a configuration of four \textit{crossing spheres} is necessary to build a crossing. Then, in order to connect all the crossings, we need to add several chains of \textit{connecting spheres}. Our aim is to determine a method for constructing necklaces where the length of the chains of spheres is as short as possible. The method developed in \cite{RR20} produces necklaces representations with one connecting sphere per crossing. This yields the upper bound $\mathrm{ball}(L)\le5\mathrm{cr}(L)$. However,  in some particular cases, we have noticed that we can modify the given packing to avoid all the connecting spheres. The method presented in this paper, based on orthocubic representations, will allow us to do so for an infinite family of alternating algebraic links, which includes the family of rational links, alternating Pretzel links or, more generally, alternating Montesinos links (see Figure \ref{fig:neck41}). Thus, we shall prove in Theorem \ref{thm:rationalball} the validity of the inequality in the Conjecture \ref{conj:4n} for this family. 

Furthermore, we are able to construct orthocubic representations that contain crossing-spheres shared by two crossings (see Figure \ref{fig:p323}). This reduce the total number of spheres to strictly fewer than 4cr($L$).  However, this phenomenon appears to disrupt the alternation between the crossings. This motivates the equality of Conjecture \ref{conj:4n}, akin to other conjectures concerning similar geometric knot invariants such as the \textit{ropelength} \cite{diao2022ropelength}, where alternation becomes a geometric constrain. We thus believe that the ball number of well-known nonalternating families, such as the torus links $T(p,q)$ with $p>q>2$, would be stritctly inferior to 4 times their crossing number. It is worthing noting that there are currently no nontrival lower-bounds on the ball number.

\subsubsection{A new visualization of the slope of rational tangles}
The correspondance between rational tangles and $\mathbb Q\cup\{\infty\}$ is well-known. Orthocubic representations provide a novel lens through which one may highlight this correspondence. Specifically, we demonstrate that the \textit{slope} of rational tangle--meaning the corresponding rational number-- can be derived from the coordinates of the intersection point of an orthocubic representation of the rational tangle and a specific circle within a cubic circle packing (see Theorem \ref{thm:slope}).

\subsubsection{Primitive solutions of a Diophantine equation}
Through the synthesis of the coordinates of the intersection point mentioned above and the inversive coordinates of spheres, we will discover an infinite set of primitive solutions for the Diophantine equation $x^4+y^4+z^4=2t^2$ (Corollary \ref{cor:diaphortho}).\\

\subsection{Organization of the paper}
The paper is structured as follows:
\medskip

 In Section \ref{sec:preliminaries}, we provide the necessary background on sphere packings and rational tangles.
 \medskip
 
In Section \ref{sec:braids}, we explore the realization of links in the tangency graph of all regular crystallographic sphere packings. We also discuss the optimality of the number of spheres used in the orthoplicial case.
  \medskip
  
In Section \ref{sec:cubicsection}, we introduce and examine orthocubic representations of rational links, and show the existence of such representations for algebraic links.
 \medskip
 
 Finally, in Section \ref{sec:points}, we delve into a geometric visualization of rational tangles and their connection to solutions of Diophantine equations.

\subsection{Acknowledgements} We express our gratitude to Alex Kontorovich for engaging discussions on Apollonian packings. We would also like to thank the referees for their valuable insights and constructive criticisms, which greatly enhanced the presentation of this work.

\section{General background}\label{sec:preliminaries} 
\subsection{Inversive coordinates}
An \textit{oriented hypersphere}, or simply \textit{sphere}, of $\wrd:=\mathbb{R}^d\cup\{\infty\}$, is the image of a spherical cap of $\mathbb S^d$ under the stereographic projection. Each sphere $S$ is uniquely defined by its center $c\in\wrd$ and its bend $b\in\mathbb R$ (the recripocal of the \textit{oriented} radius). If $S$ is a half-space, it is instead defined by its normal vector $\widehat n\in \mathbb S^{d-1}$, pointing to the interior, and the signed distance $\delta\in\ru$ between its boundary and the origin. The \textit{inversive coordinates} of $S$ are represented by the $(d+2)$-dimensional real vector
\begin{align}\label{eq:invcoord}
	\mathbf i(S)=
	\begin{cases}
		\left(bc,\dfrac{\overline b-b}{2},\dfrac{\overline b+b}{2}\right)^T&\text{ if }b\not=0,  \\
		\quad\\
		(\widehat n,\delta,\delta)^T&\text{otherwise}, \\
	\end{cases}
\end{align}
where $\overline{b}=b\|c\|^2-\frac{1}{b}$ is the \textit{co-bend} of $S$. The co-bend is the bend of $S$ after inversion through the unit sphere.  A point $P\in \widehat{\mathbb R^d}$ can be considered as a sphere with infinite bend. By taking the limit, we can define its inversive coordinates as
\begin{eqnarray}\label{eq:invcoordpoint}
	\mathbf{i}(P)&=&
	\left\lbrace\begin{array}{lll}
		\left(P,\dfrac{\|P\|^2-1}{2},\dfrac{\|P\|^2+1}{2}\right)^T&\text{ if } P\not=\infty,\\
		\quad\\
		\left(\mathbf 0_d,1,1\right)^T&\text{otherwise}.  \\
	\end{array}\right.
\end{eqnarray}
The inversive coordinates of points are \textit{homogeneous}, meaning that for every $\lambda\not=0$, $\lambda\mathbf{i}(P)$ are valid inversive coordinates of the same point \cite{wilker}. 
The \textit{inversive product} of two spheres or points $S,S'$ of $\wrd$ is the real value
\begin{align}\label{eq:invprod}
	\langle S,S'\rangle= \mathbf i(S)^T\mathbf Q_{d+2}\mathbf i(S)
\end{align}
where $\mathbf Q_{d+2}$ is the diagonal matrix $\mathrm{diag}(1,\ldots,1,-1)$ of size $d+2$. The inversive product of two spheres encodes their relative position \cite{RR20}. Furthermore, for every sphere $S$ or point $P$ of $\wrd$, we have
\begin{align}\label{eq:invprodpoint}
	\langle S,S\rangle= 1 \quad \text{and} \quad  \langle P,P\rangle= 0.
\end{align}
 An arrangement of spheres $\mathcal S$ in $\wrd$, possible infinite, is a \textit{packing} if the interiors of every two spheres are mutually disjoint. The group of Möbius transformations of $\wrd$ preserves the inversive product and acts linearly on the inversive coordinates as an orthogonal subgroup of $\mathrm {SL}_{d+2}(\mathbb R)$ with respect to $\mathbf Q_{d+2}$. In particular, the inversion\footnote{Throughout the paper, we shall use upper-case letters to denote spheres and lower-case letters to denote inversions.} $s$ through a sphere $S$  transforms the inversive coordinates through left multiplication with the matrix
\begin{align}\label{eq:invmatrix}
	\mathbf S= \mathbf{I}_{d+2}-2\mathbf i(S)^T \mathbf i(S)\mathbf{Q}_{d+2}
\end{align}
where $\mathbf{I}_{d+2}$ is the identity matrix of size $d+2$. 

\subsection{Polytopal sphere packings}
The \textit{polar}  of a subset $X\subset\mathbb R^{d+1}$  is the subset 
$X^*=\{u\in\mathbb R^{d+1}\mid \langle u, v\rangle \leq 1\text{ for all }v\in X\}$. The \textit{stereographic sphere} of a point $v\in \mathbb R^{d+1}$ \textit{outside} the unit sphere $\sd$ (i.e. with $\|v\|>1$) is the sphere $S_v$ of $\wrd$ obtained by the stereographic projection of the spherical cap $\{-v\}^*\cap \mathbb S^d$.  The \textit{arrangement projection} of polytope $\P\subset\mathbb R^{d+1}$ whose vertices are outside $\sd$ is defined as the arrangement made by the  stereographic spheres of the vertices of $\P$.
\medskip

 A polytope is termed \textit{edge-scribed} if its edges are tangent to the unit sphere, and \textit{edge-scribable} if it admits an edge-scribed realization. An edge-scribed polytope is \textit{canonical} if the barycenter of the contact points with the unit sphere is the origin.  In dimension $d\ge3$, all edge-scribed realizations of an edge-scribable $d$-polytope $\P$ are equivalent, up to Möbius transformations, to a unique canonical realization $\P_0$.     

\medskip

The arrangement projection of an edge-scribed polytope is a packing. Conversely, we say that a sphere packing $\S_\P$ in $\wrd$  with  $d\ge2$, is \textit{polytopal} if there exists an edge-scribable $(d+1)$-polytope $\P$ and a Möbius transformation $\mu$ such that $\S_\P=\mu\cdot \mathcal S_{\P_0}$, where $\mathcal S_{\P_0}$ is the arrangement projection of the canonical realization. The combinatorial structure of $\S_\P$ is encoded by the corresponding edge-scribable polytope  $\P$. The vertices and the edges of $\P$ correspond bijectively to the spheres and tangency relations of $\S_\P$. The facets of $\P$ correspond to the \textit{dual spheres} of $\S_\P$, which form the \textit{dual arrangement} $\S_\P^*:=\mu\cdot\mathcal S_{\P_0^*}$.
\medskip

 The \textit{Apollonian arrangement} of $\S_\P$ is defined as the orbit space $\mathscr{P}(\S_\P):= \langle \S_\P^*\rangle\cdot \S_\P$, where  $\langle \S_\P^*\rangle$ denotes the group generated by inversions through the dual spheres. Polytopal sphere packings and their endowed structures are unique up to Möbius transformations (see \cite{RR21_1} for more details).

\subsection{The regular crystallographic sphere packings}

In dimension two, Apollonian arrangements of polytopal circle packings are also packings, but this is not true in general \cite{RR21_1}. In higher dimensions, Apollonian arrangements that are packings belong to the family of \textit{crystallographic sphere packings} introduced by Kontorovich and Nakamura in \cite{KontorovichNakamura}. These are dense, infinite sphere packings obtained as the orbit space   $\mathscr P=\langle \widetilde {\mathcal{S}}\rangle \cdot \mathcal S$, where $\mathcal S$ is a finite sphere packing called the \textit{cluster}, $\langle \widetilde {\mathcal{S}}\rangle$ is a geometrically finite subgroup of the group of Möbius transformations generated by inversions through a finite arrangement of spheres $\widetilde {\mathcal{S}}$, called the \textit{co-cluster}. This satisfies that every sphere of $\mathcal S$ is disjoint, tangent or orthogonal to every sphere of $\widetilde {\mathcal{S}}$. Crystallographic sphere packings exist only in dimensions $2\le d \le18$ \cite{kapovich2023superintegral}. In \cite{rasskin2021regular}, the second author provided the following enumeration of all the \textit{regular} crystallographic sphere packings $\mathscr{P}_{\{p_1,\ldots,p_d\}}:= \langle \S_\P^*\rangle\cdot \S_\P$, where $\P$ is the regular polytope with Schläfli symbol $\{p_1,\ldots,p_d\}$.
\begin{enumerate}
	\item [$(d=2)$] $	\mathscr{P}_{\{3,3\}}, 	\mathscr{P}_{\{3,4\}},	\mathscr{P}_{\{4,3\}},	\mathscr{P}_{\{3,5\}},	\mathscr{P}_{\{5,3\}},	$
	\item [$(d=3)$] $	\mathscr{P}_{\{3,3,3\}}, 	\mathscr{P}_{\{3,3,4\}},	\mathscr{P}_{\{4,3,3\}}, 	\mathscr{P}_{\{3,4,3\}},	\mathscr{P}_{\{5,3,3\}},$
	\item [$(d=5)$] $	\mathscr{P}_{\{3,3,3,3,4\}}$.
\end{enumerate}

 Among the eleven regular crystallographic packings given above, two of them are of special relevance for this paper: the \textit{cubic Apollonian circle packing} $\mathscr{P}_{\{4,3\}}$ \cite{stange2015bianchi}, and the \textit{orthoplicial Apollonian sphere packing} $\mathscr{P}_{\{3,3,4\}}$ \cite{nakamura2014localglobal}. These correspond to the Apollonian arrangements of a cubic circle packing  $\C_{\{4,3\}}$ and an orthoplicial sphere packing $ \S_{\{3,3,4\}}$, respectively. 
 
 \begin{figure}[H]
 	\centering
 	
 	\begin{tikzpicture}
 		
 		\begin{scope}[xshift=-5.3cm]
 			\begin{scope}[scale=.425]
 				\draw (2.414,2.414) circle (2.414cm)  node {\Large $123$};
 				\draw (-2.414,2.414) circle (2.414cm) node {\Large$\overline{1}23$};
 				\draw (2.414,-2.414) circle (2.414cm) node {\Large$1\overline{2}3$};
 				\draw (-2.414,-2.414) circle (2.414cm) node {\Large$\overline{12}3$};
 				\draw (.414,.414) circle (.414cm) node {\tiny $12\overline{3}$};
 				\draw (-.414,.414) circle (.414cm) node {\tiny $\overline12\overline{3}$};
 				\draw (.414,-.414) circle (.414cm) node {\tiny $1\overline{23}$};
 				\draw (-.414,-.414) circle (.414cm) node {\tiny $\overline{123}$};
 			\end{scope}
 		\end{scope}
 		
 		\begin{scope}[xshift=-0cm]
 			
 			\begin{scope}[scale=.425]
 				\draw (2.414,2.414) circle (2.414cm) ;
 				\draw (-2.414,2.414) circle (2.414cm) ;
 				\draw (2.414,-2.414) circle (2.414cm) ;
 				\draw (-2.414,-2.414) circle (2.414cm) ;
 				\draw (.414,.414) circle (.414cm) ;
 				\draw (-.414,.414) circle (.414cm);
 				\draw (.414,-.414) circle (.414cm);
 				\draw (-.414,-.414) circle (.414cm);
 				
 				\draw[blue] (1.414,0) circle (1cm)  node[fill=white,inner sep=0pt]  {\Large$1$};
 				\draw[blue] (-1.414,0) circle (1cm)  node[fill=white,inner sep=0pt] {\Large$\bar1$};
 				\draw[blue] (0,1.414) circle (1cm)  node[fill=white,inner sep=0pt]  {\Large$2$};
 				\draw[blue] (0,-1.414) circle (1cm)  node[fill=white,inner sep=0pt] {\Large$\bar2$};
 				\draw[blue] (0,0) circle (2.414cm)  node at (-1.7,-1.7)[fill=white,inner sep=0pt]  {\Large$3$};
 				\draw[blue] (0,0) circle (0.414cm)  node [fill=white,inner sep=0pt] {\footnotesize$\bar3$};

 			\end{scope}
 		\end{scope}
 		\begin{scope}[xshift=5.3cm]
 			\node at (0,0) {\includegraphics[align=c,width=0.32\textwidth]{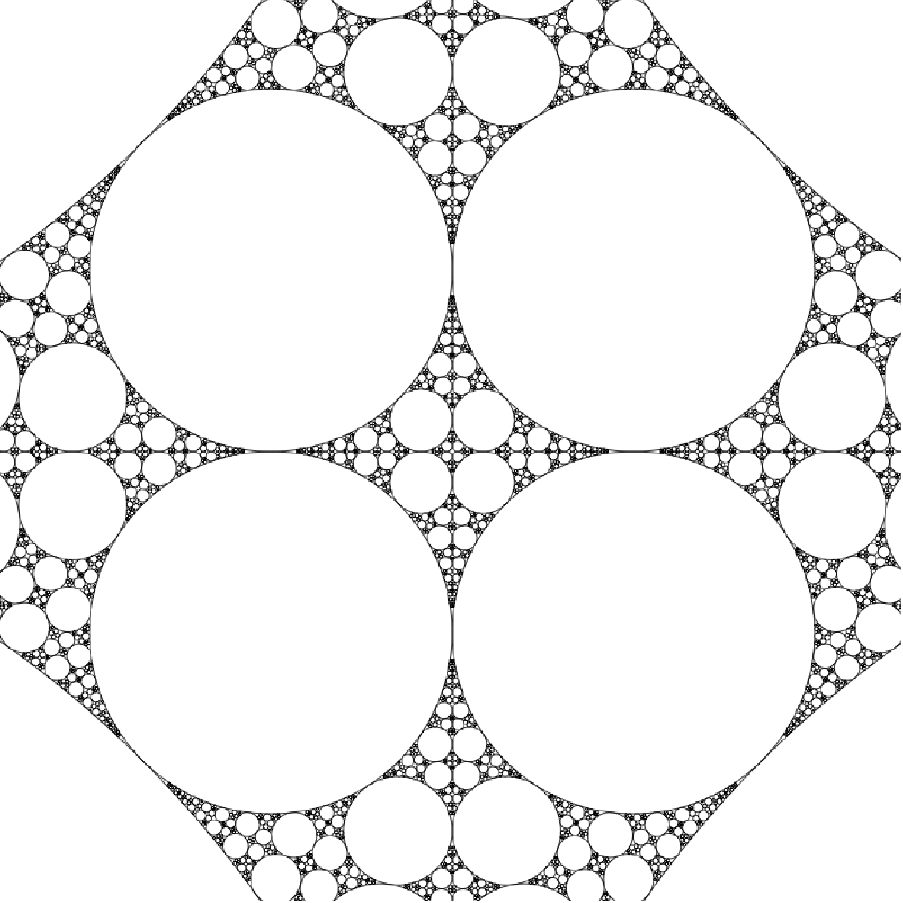}};
 		\end{scope}
 		
 	\end{tikzpicture}

 	\caption{(Left) $\C_{\{4,3\}}$ with an antipodal labelling; (center)  $\C_{\{4,3\}}$ with its dual $\C_{\{4,3\}}^*$ in blue; (right) the corresponding cubic Apollonian packing $\mathscr P_{\{4,3\}}$.}
 	\label{fig:P43}
 \end{figure}
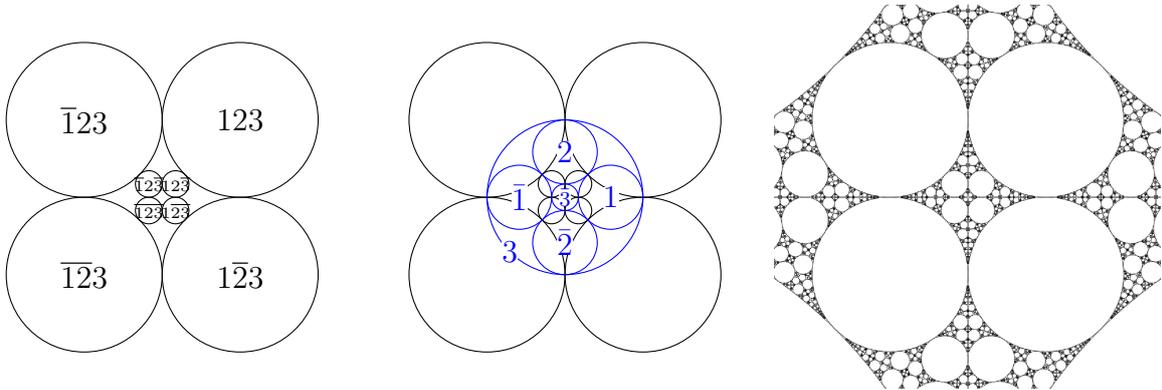
 
 It is noteworthy that the dual arrangement $\C_{\{4,3\}}^*$ is the arrangement projection of an edge-scribed octahedron, and hence a packing. However, $\S_{\{3,3,4\}}^*$ is the arrangement projection of a hypercube, which is not edge-scribed but \textit{ridge-scribed} \cite{chenpadrol}. This implies that $\S_{\{3,3,4\}}^*$ is not a packing, as the interior of its spheres overlap (see Figure \ref{fig:P43}). 
\begin{figure}[H]
	\centering
	
	\begin{tikzpicture}

			\begin{scope}[yshift=-6cm]
			\begin{scope}
				\begin{scope}[xshift=-5.cm]
					\node at (0,0) {\includegraphics[trim=80 0 0 0,clip,align=c,height=0.45\textwidth]{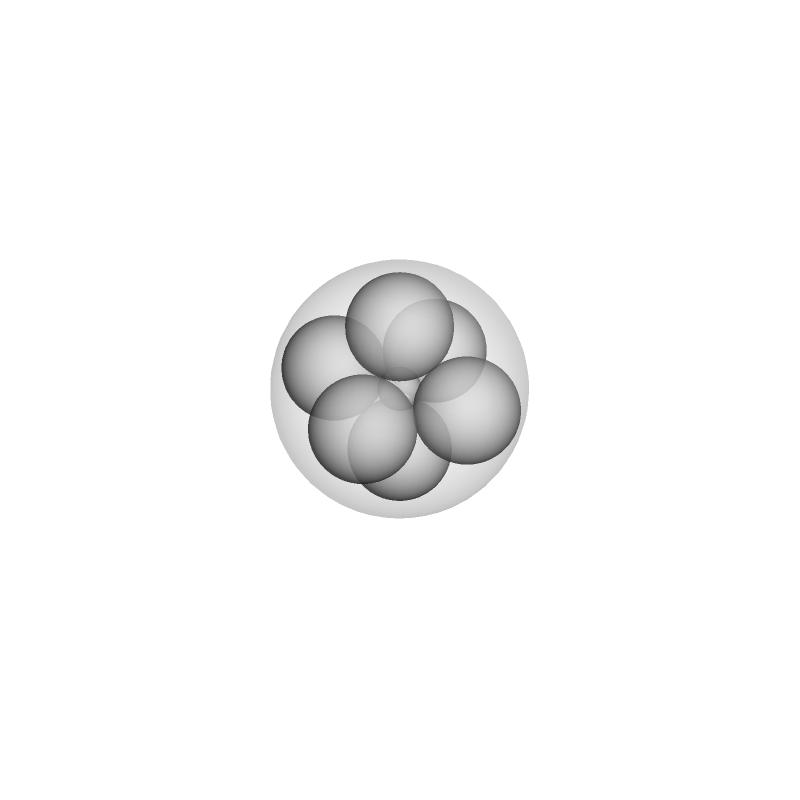}};
					
					\node at (-.15,-.13) {\Large$2$};
					\node at (-1.5,.45) {$\overline2$};
					\node at (-1.2,-.3) {\Large$1$};
					\node at (-.2,.6) {$\overline1$};
					\node at (-.8,.7) {\Large$3$};
					\node at (-.7,-.7) {$\overline3$};
					\node at (-1.6,-.9) {\Large$4$};
					\node at (-.75,.07) {\tiny$\overline 4$};
					
				\end{scope}

				\begin{scope}[xshift=5.5cm]
					\node at (0,0) {\includegraphics[align=c,height=0.45\textwidth]{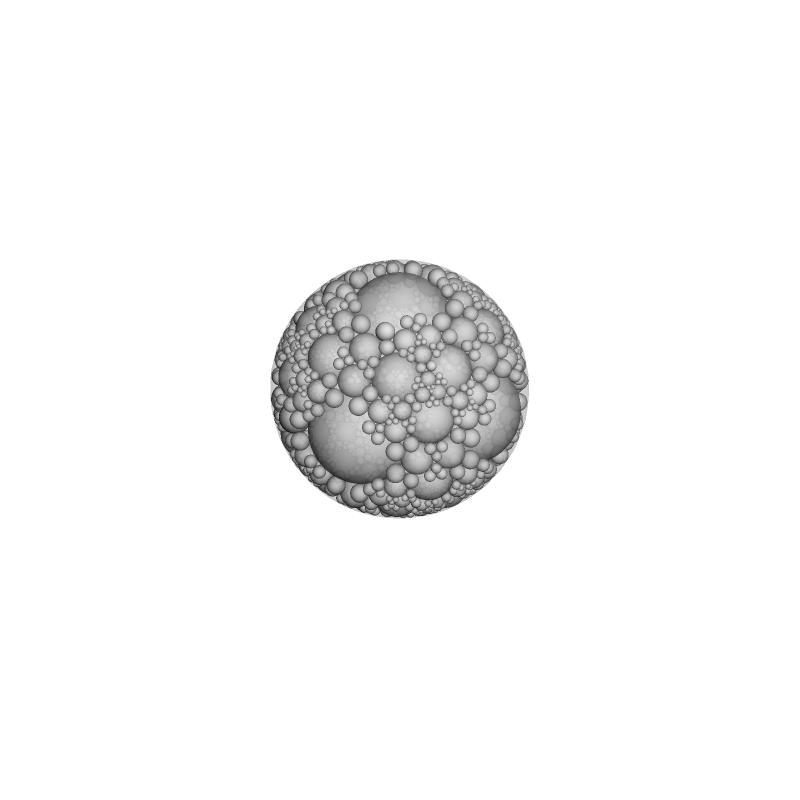}};
					
				\end{scope}
			\end{scope}
			
			\node at (0,0) {\includegraphics[align=c,height=0.45\textwidth]{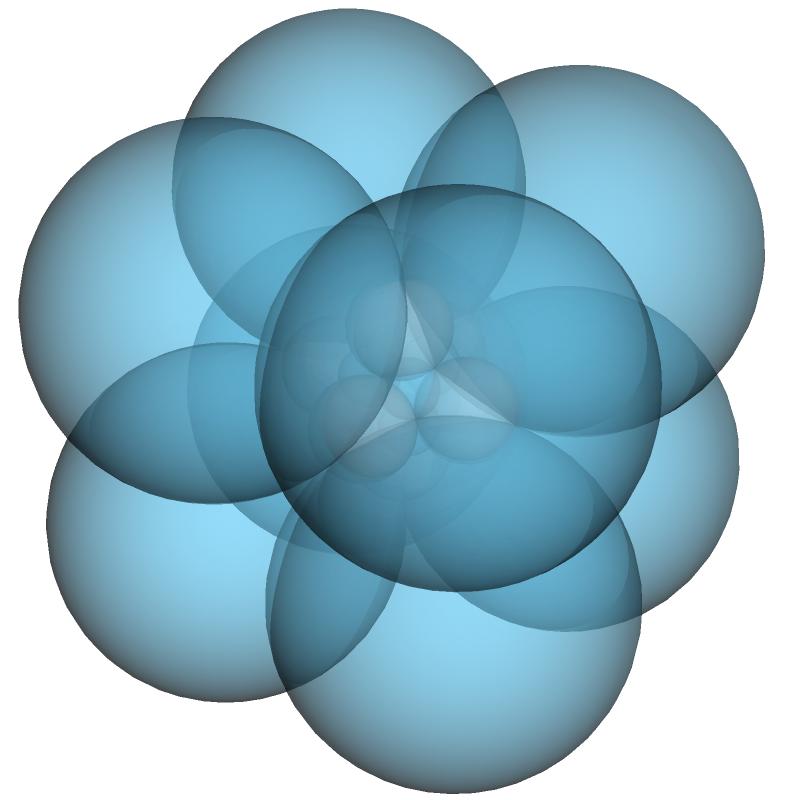}};
			\node at (.35,0) [blue]{\Large$1234$};
			\node at (.35,-2.4) [blue]{\Large$12\overline34$};
			\node at (-2.1,-1.5) [blue]{\Large$1\overline{23}4$};
			\node at (-2.1,0.8) [blue]{\Large$1\overline{2}34$};
			\node at (-.5,2.4) [blue]{\Large$\overline{12}34$};
			\node at (1.95,1.8) [blue]{\Large$\overline{1}234$};
			\node at (1.95,-.6) [blue]{\Large$\overline{1}2\overline34$};
		\end{scope}
		
	\end{tikzpicture}

	\caption{
	(Left) $\S_{\{3,3,4\}}$ with an antipodal labelling; (center)  $\S_{\{3,3,4\}}$ with its dual $\S_{\{3,3,4\}}^*$ in blue; (right) the corresponding orthoplicial Apollonian packing $\mathscr P_{\{3,3,4\}}$.
	}
	\label{fig:P334}
\end{figure}
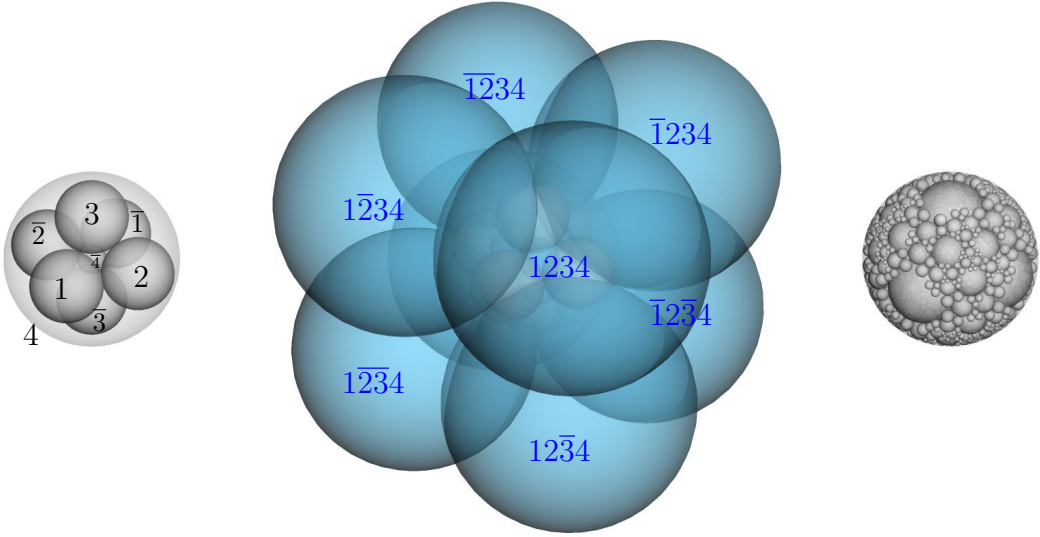

As shown in Figures \ref{fig:P43} and \ref{fig:P334}, we shall use an \textit{antipodal labelling} \cite{montejano2022self} to index the circles/spheres and the inversions of the packings $\C_{\{4,3\}}$ or $\S_{\{3,3,4\}}$ and their dual arrangements. In these labelling scheme, the labels $i$ and $\bar i:=-i$ correspond to a pair of antipodal vertices in the corresponding polytope. The vertices of the $d$-cross polytope (octahedron and orthoplex for $d=3,4$ respectively) will be labelled by $\{1,\ldots,d,\overline1,\ldots,\overline d\}$. The vertices of the $d$-cube will be labelled by concatenating the labels of the vertices incident to the corresponding facet of the dual. 
 \medskip
 
 The \textit{symmetry group} $\mathfrak{S}(\mathcal S_\P)$ of polytopal sphere packing is the stabilizer subgroup of the group of Möbius transformations for $\S_\P$. This group is isomorphic to the symmetry group of $\P$ (or its dual). The symmetry group of the $d$-cross polytope is the finite Coxeter group $B_d$ which, under the antipodal labelling, acts as the group of \textit{signed permutations} of $\{1,\ldots,d,\overline 1,\ldots,\overline d\}$. We shall denote by $r_{ij}$ (resp. $\widehat r_{ij}$)  the symmetry of the octahedron (resp. orthoplex) corresponding to the signed permutation $(ij)(\overline{ij})$.
 \medskip
 
  Different group actions might produce the same crystallographic packing. For instance, the  regular crystallographic sphere packing $\mathscr{P}_{\{p_1,\ldots,p_d\}}= \langle \S_\P^*\rangle\cdot \S_\P$ can be equally obtained as the orbit space $\Gamma_{\{p_1,\ldots,p_d\}}\cdot\{S_v\}$, where $\Gamma_{\{p_1,\ldots,p_d\}}:=\mathfrak{S}(\mathcal S_\P)\rtimes \langle \mathcal S_\P^*\rangle $ and $S_v$ is any sphere of  $\mathcal S_\P$. The group $\Gamma_{\{p_1,\ldots,p_d\}}$ is called the \textit{full symmetry group} of $\mathscr{P}_{\{p_1,\ldots,p_d\}}$ and is the stabilizer subgroup of the group of Möbius transformations for $\mathscr{P}_{\{p_1,\ldots,p_d\}}$. It is generated by a set of \textit{fundamental symmetries} generating $\mathfrak{S}(\mathcal S_\P)$, plus one inversion through a dual sphere. Then, $\mathscr{P}_{\{4,3\}}=\Gamma_{\{4,3\}}\cdot\{C_{123}\}$ and $\mathscr{P}_{\{3,3,4\}}=\Gamma_{\{3,3,4\}}\cdot\{S_1\}$ where  $\Gamma_{\{4,3\}}$ and $\Gamma_{\{3,3,4\}}$ are the following hyperbolic Coxeter groups  \cite{rasskin2021regular}.
\begin{equation}
	\begin{tabular}{cccc}
		$\Gamma_{\{4,3\}}=$\raisebox{-.4cm}{\begin{tikzpicture}[scale=.8]
				\draw[thick] (0,0) -- (3,0) 
				node[pos=.166,above,inner sep=1pt] {$4$}
				node[pos=.833,above,inner sep=1pt] {$\infty$}
				node[pos=0,below,inner sep=4pt] {$r_{3\overline3}$}
				node[pos=0.33,below,inner sep=4pt] {$r_{23}$}
				node[pos=0.66,below,inner sep=4pt] {$r_{12}$}
				node[pos=1,below,inner sep=4pt] {$s_{1}$};
				
				\foreach \x in {0,1,2,3}
				{
					\node at (\x,0) [circle,fill=black,inner sep=0pt,minimum size=.15cm]  {};
				}
				
		\end{tikzpicture} }&&&
		$\Gamma_{\{3,3,4\}}=$\raisebox{-.5cm}{	
			\begin{tikzpicture}[scale=.8]
			\draw[thick] (0,0) -- (4,0) 
			 
			node[pos=.625,above,inner sep=1pt] {$4$} 
			node[pos=.875,above,inner sep=1pt] {$4$}
			
			node[pos=0,below,inner sep=4pt] {$\widehat r_{12}$}
			node[pos=0.25,below,inner sep=4pt] {$\widehat r_{23}$}
			node[pos=0.5,below,inner sep=4pt] {$\widehat r_{34}$}
			node[pos=0.75,below,inner sep=4pt] {$	\widehat r_{4\overline4}$}
			node[pos=1.03,below,inner sep=6pt] {$s_{1234}$};
			;
			
			\foreach \x in {0,1,2,3,4}
			{
				\node at (\x,0) [circle,fill=black,inner sep=0pt,minimum size=.15cm]  {};
			}
		\end{tikzpicture} }
	\end{tabular}
\end{equation}

\subsection{Apollonian sections} Given an arrangement of spheres $\mathscr{P}:= \Gamma\cdot \S$ (not necessarily a packing),  an \textit{Apollonian section} is a subset $\mathscr{S}:= \Gamma'\cdot \S' \subset\mathscr P$ where $\Gamma'<\Gamma$ and $\S'\subset \S$. We say that $\mathscr S$ is \textit{geometric} if there is a sphere $\Sigma$, called the \textit{cutting sphere} (or plane), invariant under the action of $\Gamma'$ and intersecting all the spheres of $\mathscr S$. Two Apollonian sections $\mathscr S\subset\mathscr P $, $\mathscr S'\subset\mathscr P'$ are \textit{algebraically equivalent} if there is equivariant bijection between $\mathscr S$ and $\mathscr S'$ with respect to their group actions. 
\medskip

\begin{figure}[H]
	\centering
	
	\begin{tikzpicture}
		\begin{scope}
			\begin{scope}[xshift=-5.5cm]
				\node at (0,0) {\includegraphics[trim=0 0 30 0,clip,align=c,width=0.35\textwidth]{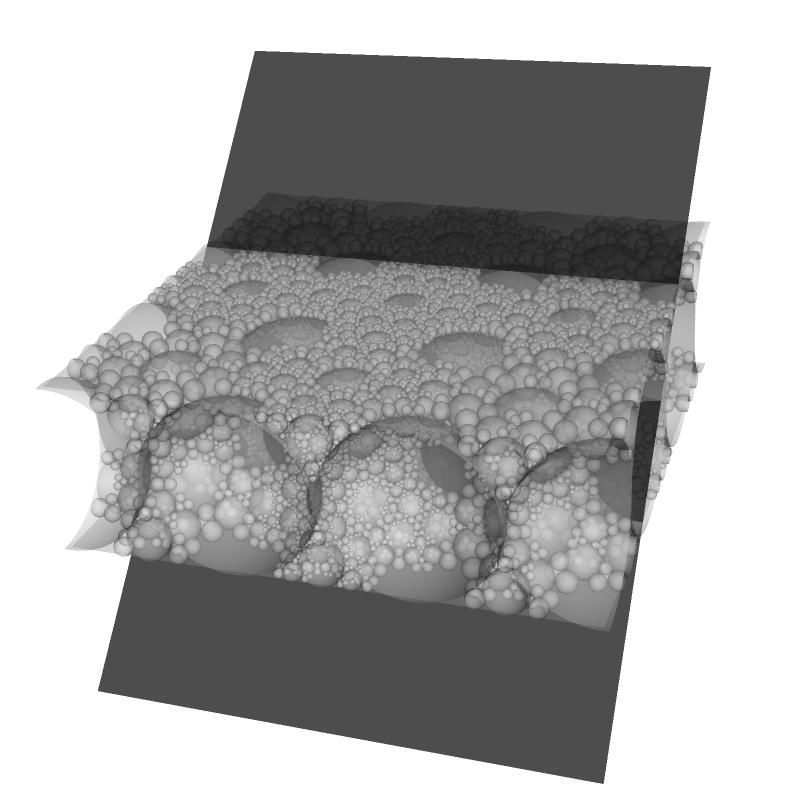}};
			\end{scope}
			\begin{scope}[xshift=-0cm]
				\node at (0,0) {\includegraphics[trim=0 0 40 0,clip,align=c,width=0.35\textwidth]{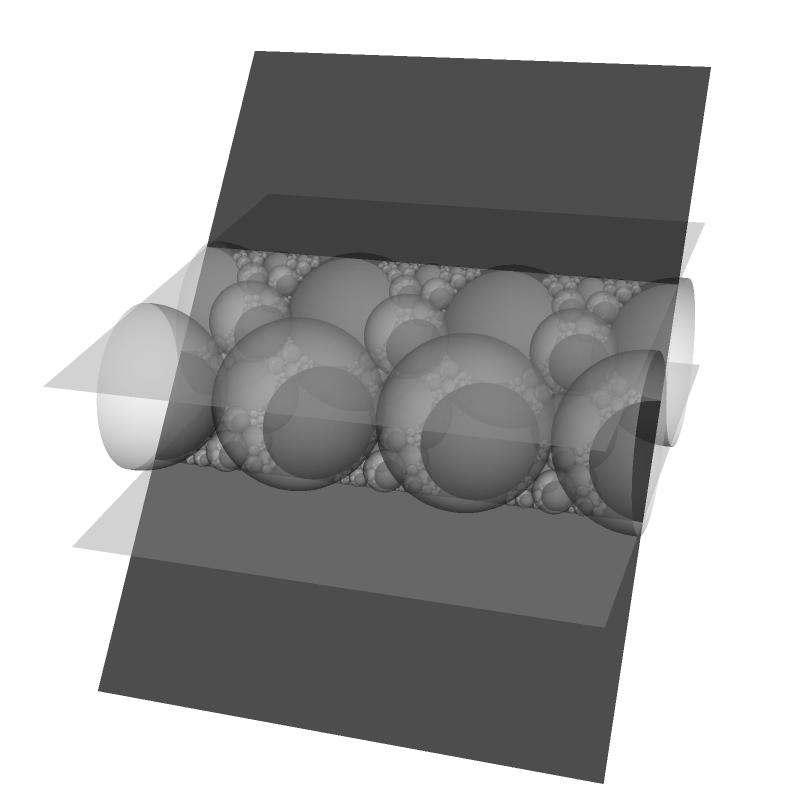}};
				
			\end{scope}
			\begin{scope}[xshift=5.cm]
				\node at (0,0) {\includegraphics[align=c,width=0.25\textwidth]{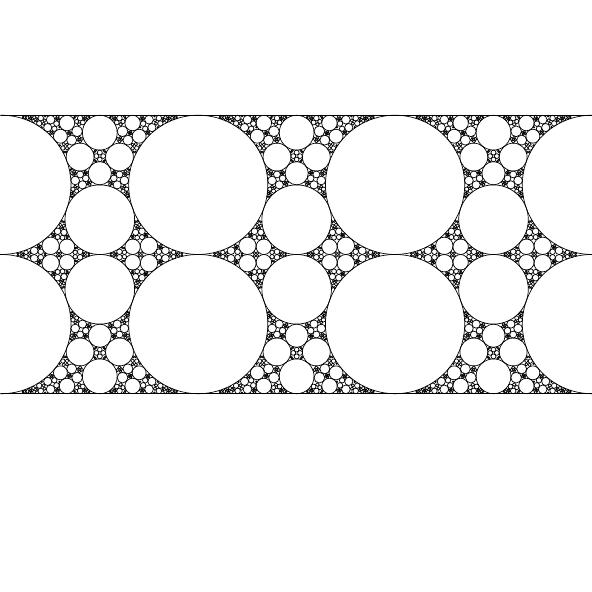}};
			\end{scope}
		\end{scope}

	\end{tikzpicture}

	\caption{(Left) $\mathscr{P}_{\{3,3,4\}}$ with the cutting plane of $\mathscr S_{\{3,3,4\}}^{\{4,3\}}$ (center), and the cubic Apollonian circle packing $\mathscr P_{\{4,3\}}$ obtained by the intersection  of $\mathscr S_{\{3,3,4\}}^{\{4,3\}}$ with the cutting plane.
	}
	\label{fig:sectionsgeo334_43}
\end{figure}

 We denote by $\mathscr S_{\{r,s,t\}}^{\{p,q\}}$ a geometric Apollonian section of $\mathscr{P}_{\{r,s,t\}}$ that is algebraically equivalent to one of the five Platonic crystallographic packings $\mathscr{P}_{\{p,q\}}$. In Figure \ref{fig:sectionsgeo334_43}, we show the cubic Apollonian section $\mathscr S_{\{3,3,4\}}^{\{4,3\}}$ described in \cite{rasskin2021regular}, which can obtained as the image of the equivariant bijection $\phi_{\{3,3,4\}}^{\{4,3\}}:\mathscr{P}_{\{4,3\}}\rightarrow\mathscr{S}_{\{3,3,4\}}^{\{4,3\}}$ induced by the morphisms
 \begin{align}\label{eq:morphisms}
 	r_{12}\mapsto \widehat{r}_{12}&&r_{23}\mapsto \widehat{r}_{23}&&r_{3\overline3}\mapsto \widehat r_{1\overline2}	\widehat r_{3\overline4}&& s_1\mapsto s_{1\overline{23}4}&& C_{123}\mapsto S_4
 \end{align}

\subsection{Algebraic links}
 \newcommand\opacb{.1}
A {\it tangle} is a pair $(S,t)$ where $S$ is a compact set of $\mathbb R^3$ homeomorphic to a $3$-ball, and $t$ is a collection $\{\gamma_1,\gamma_2,\ldots,\gamma_m\}$ of $m\ge2$ disjoint arcs contained in $S$. The arcs  $\gamma_1$ and $\gamma_2$ are open arcs whose endpoints lie on the boundary of $S$, while the rest of the arcs are closed. Two tangles $(S,t)$ and $(S',t')$ are said to be {\em equivalent} if there is an isotopy of $\mathbb R^3$ carrying $S$ to $S'$, $t$ to $t'$ and the endpoints of $(S,t)$ to the endpoints of $(S',t')$. We denote this equivalence relation by $t\simeq t'$.   Up to equivalence, we may assume that $S$ is a regular $3$-ball and the endpoints of $t$ are equally spaced along an equator $C$ of $S$. A \textit{tangle diagram} of $(S,t)$ is a regular projection of $t$ and $C$ onto the plane containing $C$, together with the crossing information. We shall denote the endpoints in a tangle diagram by the cardinal points NE, NW, SE, and SW, as shown in Figure  \ref{fig:elementarytan}. If not otherwise required, we will refer to a tangle $(S,t)$ simply as $t$.  In the following, we describe the builiding blocks for contructing algebraic tangles.

\begin{enumerate}[label=($\T_\arabic*)$]
	\setcounter{enumi}{-1}
	\item The \textit{elementary tangles} are described in Figure \ref{fig:elementarytan}.
	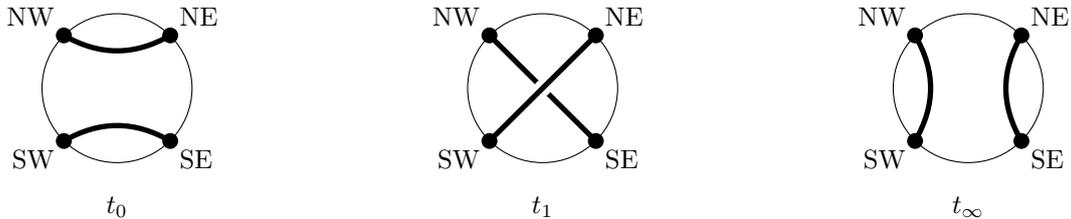
\begin{figure}[H]
		\centering
			\renewcommand{\thi}{.7mm}
\begin{tikzpicture}[scale=.7]
\begin{scope}[xshift=-8cm]
\draw  circle (1.41);

\draw[line width=\thi,bend right] (-1,1) edge (1,1);
\draw[line width=\thi,bend left] (-1,-1) edge (1,-1);

\fill[black] \NE circle (1.5mm) node[anchor= south west] {NE};
\fill[black] \NW circle (1.5mm) node[anchor= south east] {NW};
\fill[black] \SW circle (1.5mm) node[anchor=north  east] {SW};
\fill[black] \SE circle (1.5mm) node[anchor=north west] {SE};

\node at (0,-2.25) {$t_0$};
\end{scope}


\draw  circle (1.41);

\draw[line width=\thi] \NW -- \SE;
\draw[draw=white,fill=white]  circle (0.15);
\draw[line width=\thi] \SW -- \NE;

\fill[black] \NE circle (1.5mm) node[anchor= south west] {NE};
\fill[black] \NW circle (1.5mm) node[anchor= south east] {NW};
\fill[black] \SW circle (1.5mm) node[anchor=north  east] {SW};
\fill[black] \SE circle (1.5mm) node[anchor=north west] {SE};

\node at (0,-2.25) {$t_1$};

\begin{scope}[xshift=8cm]
\draw  circle (1.41);

\draw[line width=\thi,bend left] (-1,1) edge (-1,-1);
\draw[line width=\thi,bend left] (1,-1) edge (1,1);

\fill[black] \NE circle (1.5mm) node[anchor= south west] {NE};
\fill[black] \NW circle (1.5mm) node[anchor= south east] {NW};
\fill[black] \SW circle (1.5mm) node[anchor=north  east] {SW};
\fill[black] \SE circle (1.5mm) node[anchor=north west] {SE};

\node at (0,-2.25) {$t_\infty$};
\end{scope}

\end{tikzpicture}   
		\caption{The elementary tangles.}
		\label{fig:elementarytan}
	\end{figure} 
	\item The \textit{mirror} of $t$ is the tangle $-t$ given by the reflection of $t$ through the plane containing the equator.
	\item The \textit{flip} of $t$ is  the tangle $F(t)$ given by the reflection of $t$ throught the plane orthogonal to the equator and passing through the endpoints SW and NE.
	\begin{figure}[H]
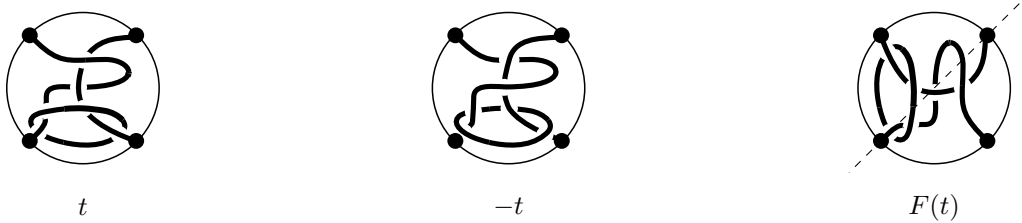

	\centering
	\includestandalone[scale=1]{tikzs/tangoperat} 
	\label{fig:tangbinop}
	\caption{The mirror and the flip of a tangle.}
\end{figure} 

	\item The \textit{addition} of two tangles $t$ and $t'$ is the tangle $t+t'$, obtained by connecting the East endpoints of $t$ to the West endpoints of $t'$.
	\item The \textit{positive (resp. negative) half-twist} of a tangle $t$ is $H^+(t):=t_1+t$ (resp. $H^{-}(t):=-t_1+t$).\\
	
		\begin{figure}[H]
		\centering
		\renewcommand{\thi}{.7mm}

\begin{tikzpicture}[scale=.7]


\begin{scope}[xshift=-8cm]

\coordinate (NE) at (1,1);
\coordinate (SE) at (1,-1);
\coordinate (NW) at (-1,1);
\coordinate (SW) at (-1,-1);

\draw [line width=\thi] circle (1.41);

\node at (0,-2.25) {$t+t'$};

\begin{scope}[xshift=-2cm]

\draw[fill=white]  circle (1.41);

\fill[black] \NW circle (1.5mm);
\fill[black] \SW circle (1.5mm);

\node at (0,0) {\large $t$};
\end{scope}


\begin{scope}[xshift=2cm]
	
	\draw[fill=white]  circle (1.41);
	
	\fill[black] \NE circle (1.5mm);
	\fill[black] \SE circle (1.5mm);
	
	\node at (0,0) {\large $t'$};
	
\end{scope}
\end{scope}


\begin{scope}
\coordinate (NE) at (1,1);
\coordinate (SE) at (1,-1);
\coordinate (NW) at (-1,1);
\coordinate (SW) at (-1,-1);

\node at (0,-2.25) {$H^+(t)$};

\draw [line width=\thi] circle (1.41);


\begin{scope}[xshift=-2cm]

	\draw[fill=white]  circle (1.41);
	
	\fill[black] \NW circle (1.5mm);
	\fill[black] \SW circle (1.5mm);
	
	\draw[line width=\thi] \NW -- \SE;
	\draw[line width=3*\thi,white] (-.6,-.6) --  (.6,.6);
	\draw[line width=\thi] \SW -- \NE;
	
\end{scope}	
	
	\begin{scope}[xshift=2cm]
		
		\coordinate (NW) at (-1,1);
		\coordinate (SW) at (-1,-1);
		\draw[fill=white] circle (1.41);
		\node at (0,0) {\large $t$};
		
		\fill[black] \NE circle (1.5mm);
		\fill[black] \SE circle (1.5mm);
	\end{scope}

\end{scope}

\begin{scope}[xshift=8cm]
	\coordinate (NE) at (1,1);
	\coordinate (SE) at (1,-1);
	\coordinate (NW) at (-1,1);
	\coordinate (SW) at (-1,-1);
	
	\node at (0,-2.25) {$H^-(t)$};
	
	\draw [line width=\thi] circle (1.41);
	
	
	\begin{scope}[xshift=-2cm]
		
		\draw[fill=white]  circle (1.41);
		
		\draw[line width=\thi] \SW -- \NE;
		\draw[line width=3*\thi,white] (-.6,.6) --  (.6,-.6);
		\draw[line width=\thi] \NW -- \SE;
		
		\fill[black] \NW circle (1.5mm);
		\fill[black] \SW circle (1.5mm);

	\end{scope}	
	
	\begin{scope}[xshift=2cm]
		
		\coordinate (NW) at (-1,1);
		\coordinate (SW) at (-1,-1);
		\draw[fill=white] circle (1.41);
		\node at (0,0) {\large $t$};
		
		\fill[black] \NE circle (1.5mm);
		\fill[black] \SE circle (1.5mm);
	\end{scope}
	
\end{scope}
\end{tikzpicture}  
		\caption{The addition and the half-twists of tangles.}
		\label{fig:tangsum}
	\end{figure}
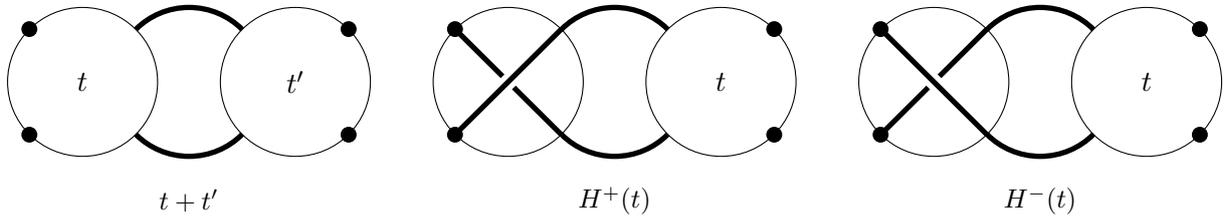 
	\item 		The \textit{closure} of a tangle $(S,t)$ is the link formed by joining the endpoints by two disjoint and unlinked paths at the exterior of $S$. Up to isotopy, there are two possible closures, the {\em numerator}  $N(t)$, obtained by joining the northern and the southern endpoints separately, and the {\em denominator} $D(t)$, obtained by joining the western and the eastern endpoints (see Figure \ref{fig:tangcloses}). 
	\begin{figure}[H]
		\centering
		\renewcommand{\thi}{.7mm}
\begin{tikzpicture}[scale=.7]
\draw  circle (1.41);
\draw[->] (1.8,0)--(3.7,0) node [below,pos=.5] {$D$};
\draw[->] (-1.8,0)--(-3.7,0) node [below,pos=.5] {$N$};
\node at (0,0) {\Large $t$};
\fill[black] \NW circle (1.5mm);
\fill[black] \SW circle (1.5mm);
\fill[black] \NE circle (1.5mm);
\fill[black] \SE circle (1.5mm);

\begin{scope}[xshift=-6cm]
\draw  circle (1.41);
\draw[line width=\thi] (-1,1) .. controls (-2,2) and (2,2) .. (1,1);
\draw[line width=\thi] (-1,-1) .. controls (-2,-2) and (2,-2) ..  (1,-1);
\node at (0,0) {\Large $t$};
\end{scope}

\begin{scope}[xshift=6cm]
\draw  circle (1.41);
\draw[line width=\thi] (-1,1) .. controls (-2,2) and (-2,-2) .. (-1,-1);
\draw[line width=\thi] (1,1) .. controls (2,2) and (2,-2) ..  (1,-1);
\node at (0,0) {\Large $t$};
\end{scope}

\end{tikzpicture}   
		\caption{The tangle closures.}
		\label{fig:tangcloses}
	\end{figure}
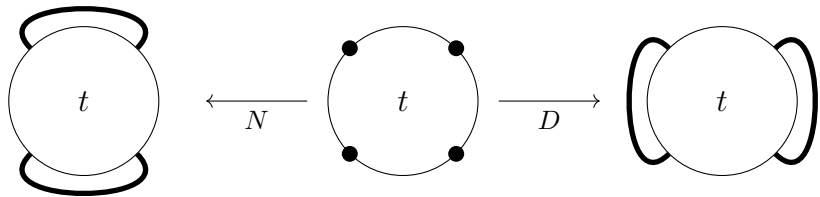 

\end{enumerate}

\begin{figure}[H]
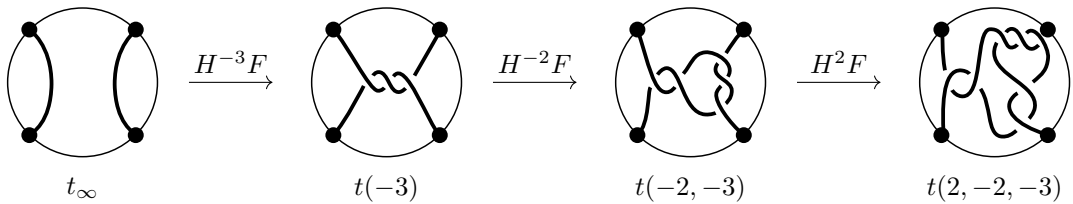

		\centering
		\includestandalone[scale=1]{tikzs/tang223all} 
		\vspace{.5cm}
		\caption{The rational tangle $t(2,-2,-3)$ obtained by the Conway's algorithm.}\label{fig:tang223}
	\end{figure}

\begin{figure}[H]
			\centering
			\includegraphics[width=.5\textwidth]{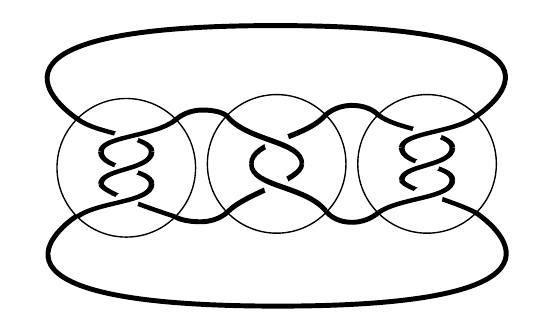}  
			\caption{The Pretzel knot $P(3,-2,3)$ which corresponds to the knot $8_{19}$ in the Alexander-Briggs notation.}\label{fig:pretzelknot}
	\end{figure}

Rational tangles were introduced by Conway in his work on enumerating and classifying knots and links \cite{conway1970enumeration}. For a given sequence of integers $a_1,\dots ,a_n$ all non-zero except maybe $a_1$, we denote by  $t(a_1,\cdots,a_n)$ the \textit{rational tangle}  given by the following Conway's algorithm \cite{cromwell2004knots} (see Figure \ref{fig:tang223}).
	\begin{align}
		t(a_1,\cdots,a_n):=H^{a_1}F\cdots H^{a_n}F(t_\infty).
	\end{align}
	The \textit{slope} of a rational tangle $t(a_1,\ldots,a_n)$ is the rational number $p/q$ obtained by the
	continued fraction expansion
	\begin{align}
		[a_1,\ldots,a_n]:=a_1+\frac1{{\ddots+\frac{1}{a_n}}}=\frac{p}{q}.
	\end{align}
	The term \textit{rational tangle} originates from the connection established by Conway \cite{conway1970enumeration}, which relates the family of tangles produced by Conways's algorithm to rational numbers. Conway's theroem states that two rational tangles are equivalent if and only if they have the same slope. We denote by $t_{p/q}$ the class of rational tangles with slope $p/q$ up to isotopy.	A \textit{rational link} is the closure of a rational tangle. \textit{Algebraic tangles} are those obtained by additions and flips of rational tangles \cite{cromwell2004knots}. Equivalently, links obtained by the closure of algebraic tangles are termed \textit{algebraic} or \textit{arborescent} \cite{gabai1986genera}. Pretzel links $P(q_1,\ldots,q_n):=N(t_{1/q_1}+\cdots+t_{1/q_n})$ are a particular case of algebraic links, as shown in Figure \ref{fig:pretzelknot}.

\section{Necklace representations in the regular crystallographic sphere packings} \label{sec:braids}
In this section, we investigate the following question: given a link $L$ and a three-dimensional crystallographic sphere packing $\mathscr P$, can we find a necklace representation of $L$ contained in $\mathscr P$? We answer this question affirmatively in the regular case. 
\medskip

First, let us introduce a concept needed for the proof. For any regular crystallographic sphere packing $\mathscr{P}_{\{p,q,r\}}$ and for any face $f$ of the corresponding regular $4$-polytope, we define the $f$-\textit{section} of $\mathscr{P}_{\{p,q,r\}}$ as the Apollonian section $\Gamma_{\{p,q,r\}}^{f}\cdot \mathcal S^f_{\{p,q,r\}}\subset\mathscr{P}_{\{p,q,r\}}$. Here, $\Gamma_{\{p,q,r\}}^{f}$ is the stabilizer subgroup of $\Gamma_{\{p,q,r\}}$ for $\{S_v\}_{v\in f}$, and  $\mathcal S_{\{p,q,r\}}^f:=\mathcal S_{\{p,q,r\}}\setminus\{S_v\}_{v\in f}$.
\begin{thm}\label{lem:ortholinks}
	Every link admits a necklace representation in all three-dimensional regular crystallographic sphere packings.
\end{thm}  

\begin{figure}[H]
	\centering
	\includegraphics[width=.4\textwidth,align=c]{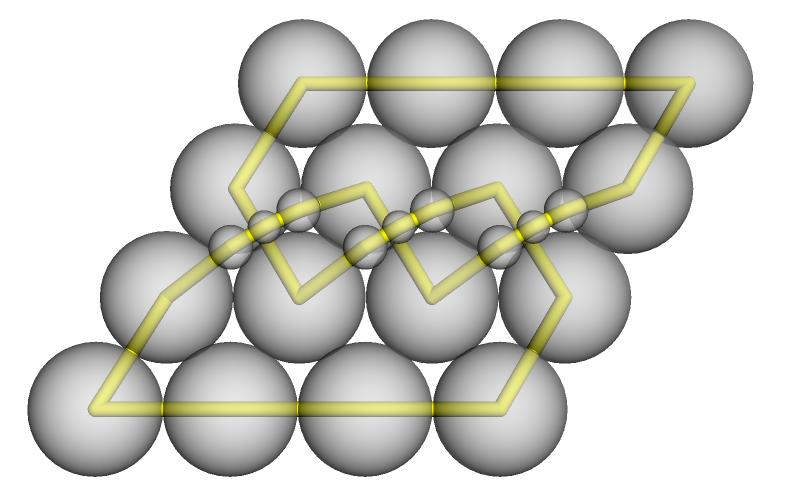} 	\includegraphics[width=.4\textwidth,align=c]{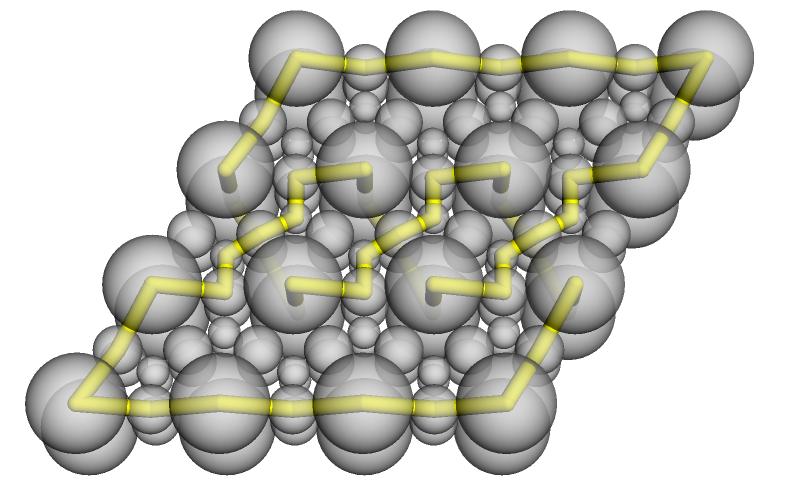} 
	
	\caption{Two necklace representations of the trefoil knot: one contained in a vertex section of $\mathscr P_{\{3,3,3\}}$ (left) and  the other in an edge section of $\mathscr P_{\{4,3,3\}}$ (right).}
	\label{fig:traing}
\end{figure}

\begin{proof} Let $L$ be a link. We start by constructing a necklace representation of $L$ in the orthoplicial Apollonian packing $\mathscr{P}_{\{3,3,4\}}$, which contains the strip orthoplicial packing $\mathcal{S}_{\{3,3,4\}}$ obtained by applying an inversion through a sphere centered at the contact point of $S_1$ and $S_2$. We consider the edge section $\Gamma_{\{3,3,4\}}^{12}\cdot \mathcal S^{12}_{\{3,3,4\}}\subset \mathscr{P}_{\{3,3,4\}}$ for the edge $(1,2)$. We have that  $\Gamma_{\{3,3,4\}}^{12}=\langle \widehat r_{34},\widehat r_{4\overline4},s_{1234}\rangle<\Gamma_{\{3,3,4\}}$, and $\mathcal S_{\{3,3,4\}}^{12}=\{S_3,S_4,S_{\overline1},S_{\overline2},S_{\overline3},S_{\overline 4}\}$. The orbits of $\{S_3,S_4,S_{\overline3},S_{\overline 4}\}$ form an infinite square grid in the tangency graph of the edge section, 
 while the orbits of $\{S_{\overline1},S_{\overline2}\}$ add addional vertices above and below each square (see Figure \ref{fig:sqpacking2}).

\begin{figure}[H]
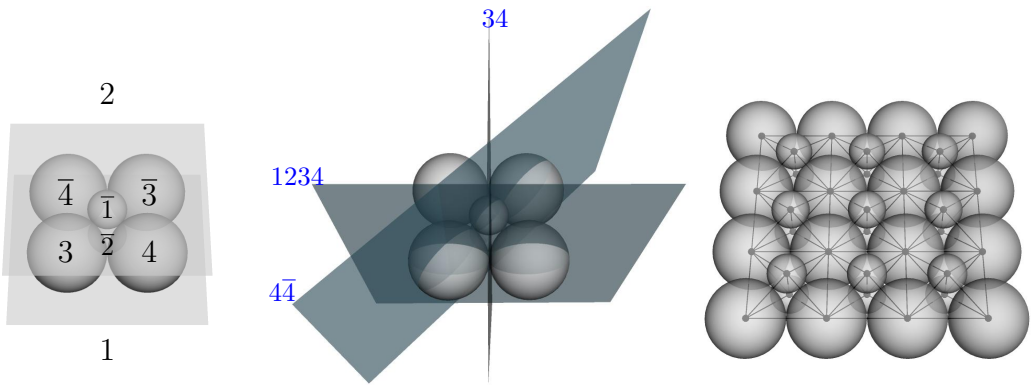

	\centering
	\includestandalone[scale=1,align=c]{tikzs/sqsection_v3}	
	\caption{(Left) The strip orthoplicial sphere packing $\mathcal S_{\{3,3,4\}}$; (center) the subset $\mathcal S_{\{3,3,4\}}^{12}$ with the walls of the generators of $\Gamma_{\{3,3,4\}}^{12}$; (right) the edge section of $\mathscr{P}_{\{3,3,4\}}$ with its tangency graph.}
	\label{fig:sqpacking2}
\end{figure}

By the well-known Alexander's Theorem \cite{alexander1923lemma},  for every link $L$ there exists a braid $\gamma$ such that its closure is isotopically equivalent to $L$. We can always draw a diagram of $\gamma$ in a regular square grid, where the crossings are placed at the intersections of the diagonals of the squares, and the remaining arcs follow the edges of the grid. This \textit{braid-grid diagram} induces a polygonal closed path in the tangency graph of the edge section, providing us with a necklace representation of $L$ in $\mathscr P_{\{3,3,4\}}$, as depicted in Figure \ref{fig:sqgrid31}. 

\begin{figure}[H]
	\centering
	\includegraphics[width=.25\textwidth,align=c]{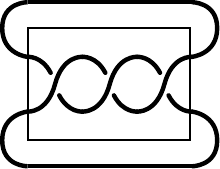}\hspace{1cm}  
\begin{tikzpicture}[scale=1.25]
\begin{scope}[xshift=-.5cm,yshift=-1.5cm,rotate=90]
\draw[help lines,step=1,dashed] (-1,0) grid (2,3); 

\begin{scope}[yshift=0cm]
\draw[line width=\thi,line cap=round] (0,0) -- (1,1);
\draw[draw,white, fill=white] (.5,.5) circle (1.2mm);
\draw[line width=\thi,line cap=round] (1,0) -- (0,1);
\draw[line width=\thi,line cap=round] (2,0) -- (2,1);
\draw[line width=\thi,line cap=round] (-1,0) -- (-1,1);
\end{scope}

\begin{scope}[yshift=1cm]
\draw[line width=\thi,line cap=round] (0,0) -- (1,1);
\draw[draw,white, fill=white] (.5,.5) circle (1.2mm);
\draw[line width=\thi,line cap=round] (1,0) -- (0,1);
\draw[line width=\thi,line cap=round] (2,0) -- (2,1);
\draw[line width=\thi,line cap=round] (-1,0) -- (-1,1);
\end{scope}

\begin{scope}[yshift=2cm]
\draw[line width=\thi,line cap=round] (0,0) -- (1,1);
\draw[draw,white, fill=white] (.5,.5) circle (1.2mm);
\draw[line width=\thi,line cap=round] (1,0) -- (0,1);

\draw[line width=\thi,line cap=round] (2,0) -- (2,1);
\draw[line width=\thi,line cap=round] (-1,0) -- (-1,1);
\end{scope}

\begin{scope}[yshift=0cm]
\draw[line width=\thi,line cap=round]  (1,0) -- (2,0);
\draw[line width=\thi,line cap=round]  (-1,0) -- (0,0);
\end{scope}

\begin{scope}[yshift=3cm]
\draw[line width=\thi,line cap=round] (1,0) -- (2,0);
\draw[line width=\thi,line cap=round] (-1,0) -- (0,0);
\end{scope}

\end{scope}

\end{tikzpicture}  \hspace{1cm}\includegraphics[height=4cm,align=c]{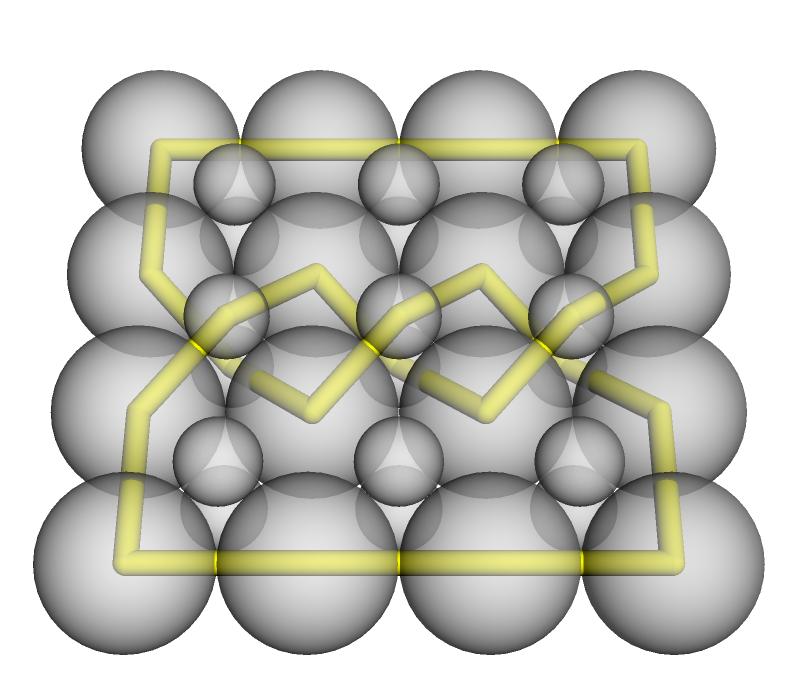} 
	\caption{(Left) A diagram of the trefoil obtained as the closure of a braid; (center) a braid-grid diagram of the same closed braid; (right) the corresponding necklace representation of the trefoil in the edge section of  $\mathscr P_{\{3,3,4\}}$.}
	\label{fig:sqgrid31}
\end{figure}
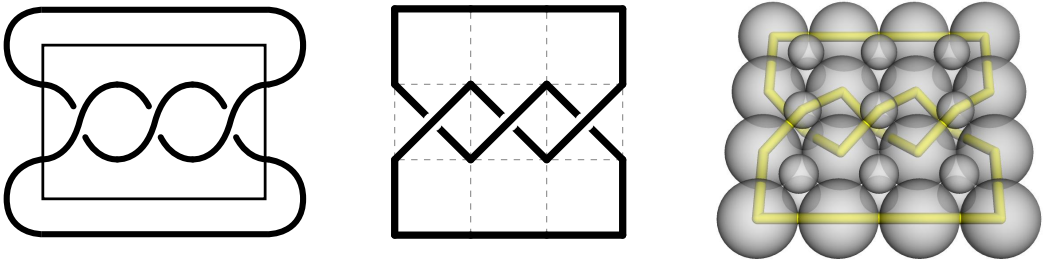
We apply a similar strategy to construct a necklace representation of any link contained in the edge section of $\mathscr{P}_{\{4,3,3\}}$, $\mathscr{P}_{\{3,4,3\}}$ and $\mathscr{P}_{\{5,3,3\}}$. The tangency graph of these three sections contains a triangular grid instead of a square grid. In this triangular grid, two tangent triangles forming a rhombus can serve the same function as a square in the braid-grid diagrams. Furthermore, for each of the three cases and each type of crossing, we can find two disjoint chains of spheres connecting the opposite vertices of the rhombus to construct the crossing. Finally, for the packing $\mathscr{P}_{\{3,3,3\}}$, the strategy is almost the same as in the previous three cases, but we must use a vertex section instead of the edge section to contruct the chains of spheres needed for the crossings (see Figure \ref{fig:traing}). 
\end{proof} 

 We observe that Theorem \ref{lem:ortholinks} can be proved without invoking Alexander's Theorem by using \textit{grid diagrams} (a slightly different notion of the braid-grid diagrams described in the proof, where the crossings are also projected onto the grid). Indeed, Cromwell proved in \cite{cromwell1995embedding} that every link admits a grid diagram. However, the necklaces induced by these diagrams are generally not optimal in terms of the number of spheres. The reason is that the chains of connecting spheres between crossings that are far apart increase the total number of spheres to a quantity greater than the upper bound given in \cite{RR20}. A similar phenomenon occurs for the necklaces induced by the closure of braid-grid diagrams of braids with more than 3 strands. However, in the particular case of $2$-braid links, we can reduce the number of spheres by modifying the necklace representation contained in the edge-section of $\mathscr{P}_{\{3,3,4\}}$.
 
\begin{figure}[H]
 	\centering

 	\includegraphics[width=.38\textwidth,align=c]{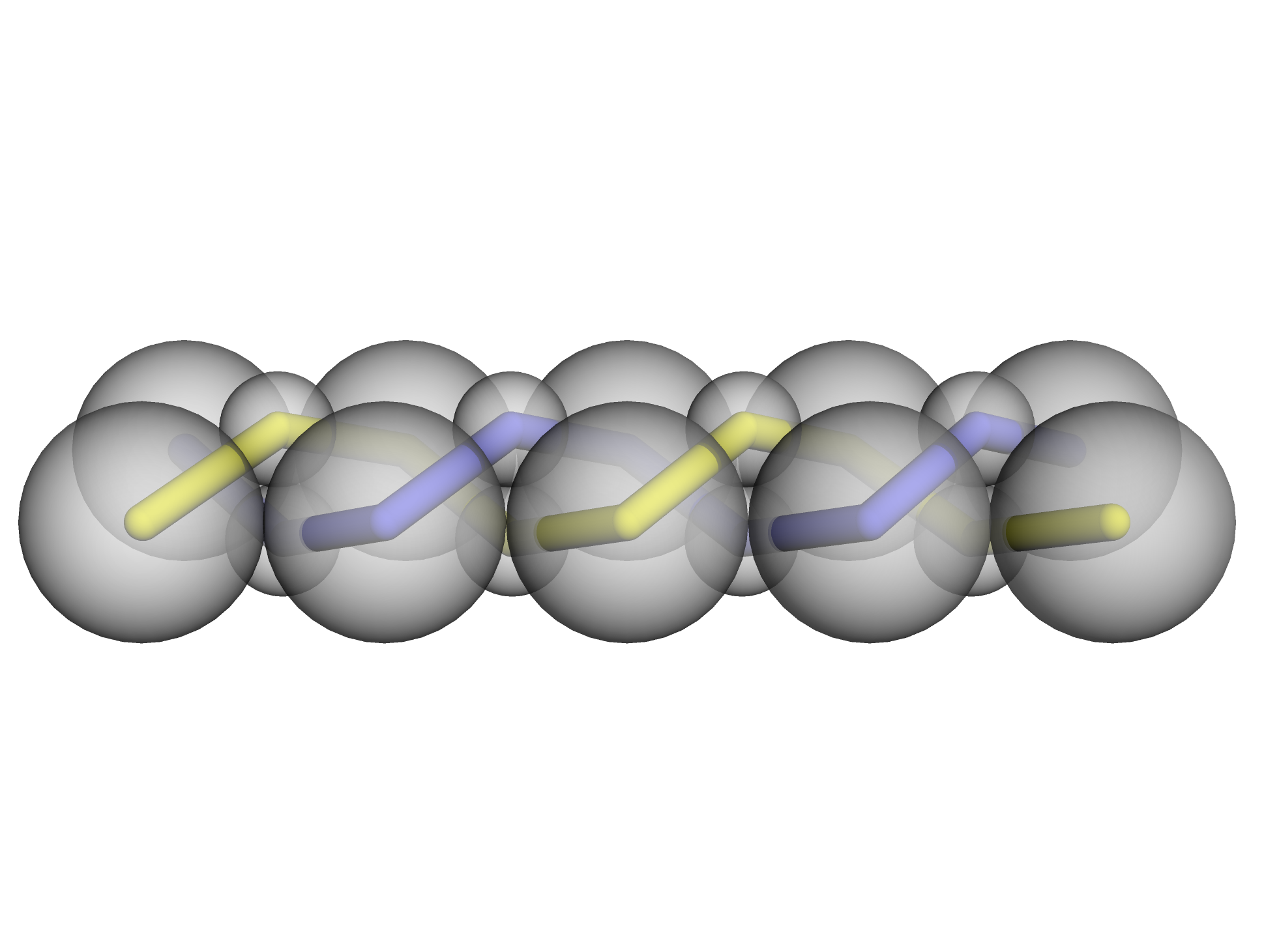} 
 	\includegraphics[width=.36\textwidth,align=c]{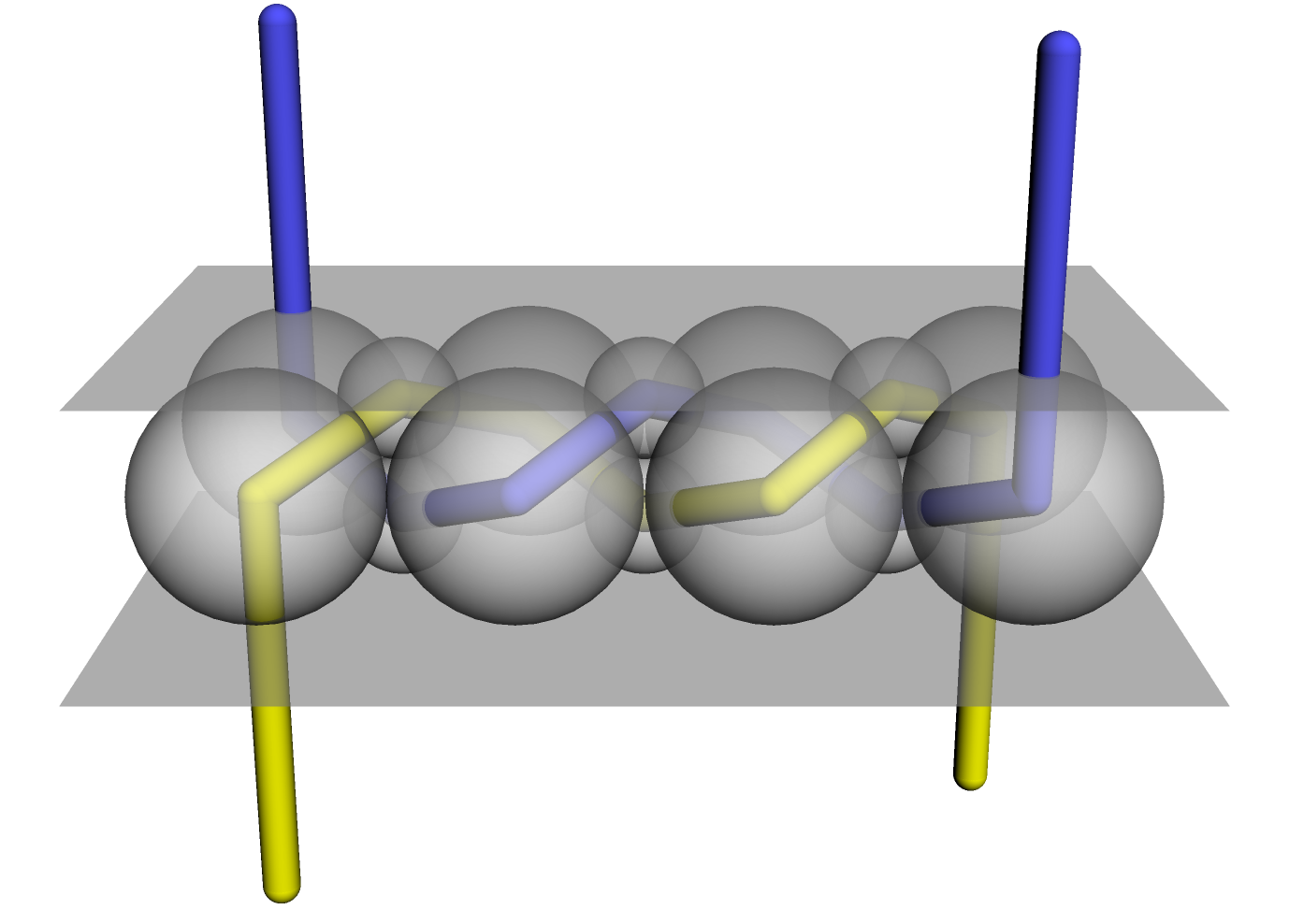} 
 	 \vspace{-.5cm}
 	\caption{(Left) A necklace representation of the alternating $2$-braid link of four crossings in the braid-grid section of $\mathscr P_{\{3,3,4\}}$; (right) the induced necklace with 16 spheres in  $\mathscr P_{\{3,3,4\}}$ obtained by using the half-spaces for the closure.}
 	\label{fig:2braid}
 \end{figure}

 \begin{cor}\label{cor:2braid} For any $2$-braid link $L$, $\mathrm{ball}(L)\le 4\mathrm cr(L).$
 \end{cor}
 \begin{proof}
 	The necklace representation in the edge section of $\mathscr{P}_{\{3,3,4\}}$ induced by the braid-grid diagram of an alternating $2$-braid with $n$ crossings consists of $4n+2$ spheres. We can then construct the closure by replacing the last 4 spheres with the two half-spaces of $\mathscr{P}_{\{3,3,4\}}$ (see Figure \ref{fig:2braid}). Since the link is alternating, $n=\mathrm{cr}(L)$.
 \end{proof}
 \vspace{-.3cm}
 
 Among the necklace representations contained in the regular crystallographic sphere packings, those produced in the orthoplicial case appear to be more optimal in terms of the number of spheres. However, we believe that the other cases could be interesting for other purposes, such as constructing $4$-polytopes containing a given link in their graph \cite{eppstein2014links} or for constructing necklace representations with horoballs as in \cite{gabai2021hyperbolic}.

\section{The orthocubic representations of algebraic links}\label{sec:cubicsection}
In Theorem \ref{lem:ortholinks}, we proved that every link is contained in the tangency graph of the three-dimensional regular crystallographic sphere packings, including the orthoplicial Apollonian packing $\mathscr{P}_{\{3,3,4\}}$. In this section, we aim to refine this result by demonstrating that the family of algebraic links is contained within the cubic Apollonian section $\mathscr{S}_{\{3,3,4\}}^{\{4,3\}}\subset\mathscr{P}_{\{3,3,4\}}$. The construction required for the proof enables us to establish the inequality of Conjecture \ref{conj:4n} for an infinite family of alternating algebraic links.

\begin{figure}[H]
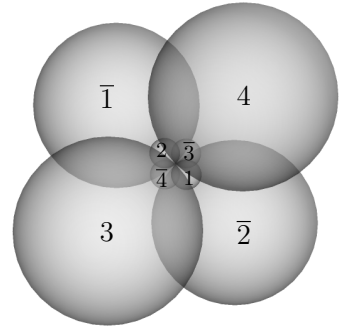

	\centering
	
	{\small
		\setlength{\tabcolsep}{2pt}
		\begin{tabular}{c|c|rrrr|rrrlrr}
			Sphere	& Bend     & \multicolumn{4}{c|}{Center }                 & \multicolumn{6}{c}{Inversive coordinates}                          \\ \hline
			$S_1$    & 	$1+1/\sqrt2$   & $(-1+\sqrt2)$ & $(\phantom-1$ & $-1$ & $1)$  & $1/\sqrt2$ & $(\phantom-1$ & $-1$ & $\phantom-1$ & $-1$ & $\sqrt2)$ \\
			$S_2$    & $1+1/\sqrt2$   & $(-1+\sqrt2)$ & $(-1$         & $1$  & $1)$  & $1/\sqrt2$ & $(-1$         & $1$  & $\phantom-1$ & $-1$ & $\sqrt2)$ \\
			$S_3$    & $1-1/\sqrt2$   & $(\phantom-1+\sqrt2)$ & $(-1$         & $-1$ & $1)$  & $1/\sqrt2$ & $(-1$         & $-1$ & $\phantom-1$ & $1$  & $\sqrt2)$ \\
			$S_4$    & $1-1/\sqrt{2}$ & $(\phantom-1+\sqrt	2)$ & $(\phantom-1$ & $1$  & $1)$  & $1/\sqrt2$ & $(\phantom-1$ & $1$  & $\phantom-1$ & $1$  & $\sqrt2)$ \\
			$S_{\overline1}$ & 	$1-1/\sqrt2$   & $(\phantom-1+\sqrt2)$ & $(-1$         & $1$  & $-1)$ & $1/\sqrt2$ & $(-1$         & $1$  & $-1$         & $1$  & $\sqrt2)$ \\
			$S_{\overline2}$ & $1-1/\sqrt2$   & $(\phantom-1+\sqrt2)$ & $(\phantom-1$ & $-1$ & $-1)$ & $1/\sqrt2$ & $(\phantom-1$ & $-1$ & $-1$         & $1$  & $\sqrt2)$ \\
			$S_{\overline3}$ & $1+1/\sqrt2$   & $(-1+\sqrt2)$ & $(\phantom-1$ & $1$  & $-1)$ & $1/\sqrt2$ & $(\phantom-1$ & $1$  & $-1$         & $-1$ & $\sqrt2)$ \\
			$S_{\overline4}$ &  $1+1/\sqrt2$   & $(-1+\sqrt2)$ & $(-1$         & $-1$ & $-1)$ & $1/\sqrt2$ & $(-1$         & $-1$ & $-1$         & $-1$ & $\sqrt2)$
		\end{tabular}
	}\includestandalone[trim=0 0 0 0,clip,align=c,scale=.9]{tikzs/orthoB1z} 
	\caption{The $z$-alternating orthoplicial sphere packing.}
	\label{fig:hoctFC}
\end{figure}

\medskip
In the following, we will establish the realization of the polytopal packings by defining $\mathcal{C}_{\{4,3\}}$ as the cubic circle packing shown in Figure \ref{fig:P43}, and $\mathcal{S}_{\{3,3,4\}}$ as the $z$-\textit{alternating} orthoplicial sphere packing described in Figure \ref{fig:hoctFC}. Notably, in the latter packing, a sphere sphere $S_i$
has positive label if it is centered above the plane $\{z=0\}$, and it has negative label if it is centered below. The cubic $\mathscr{P}_{\{4,3\}}$ and the orthoplicial $\mathscr{P}_{\{3,3,4\}}$ crystallographic packings will be those spanned by $\mathcal{C}_{\{4,3\}}$ and $\mathcal{S}_{\{3,3,4\}}$, respectively.

\begin{figure}[H]
	\centering
	\includegraphics[width=.38\textwidth,align=c]{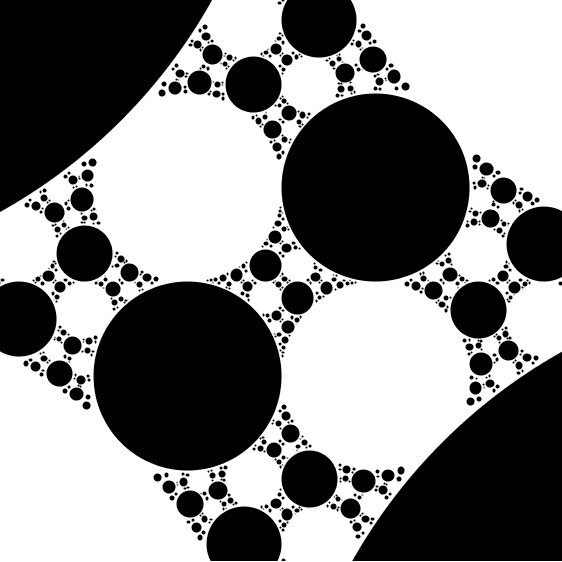} 
	\includegraphics[width=.5\textwidth,align=c]{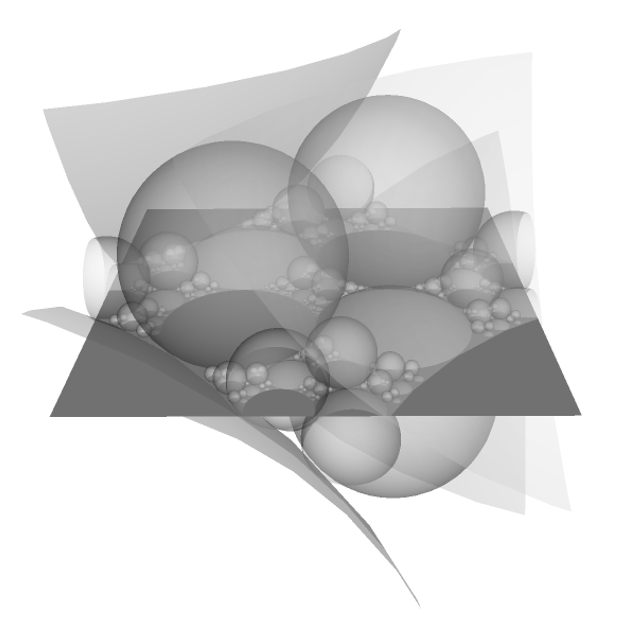} 	
	\caption{(Left) $\mathscr{P}_{\{4,3\}}$ with the $z$-coloring, (right)  $\mathscr{S}_{\{3,3,4\}}^{\{4,3\}}$ with its cutting plane $\Sigma$.}\label{fig:BWs}
\end{figure}

The geometric framework of this section will be the cubic Apollonian section $\mathscr{S}_{\{3,3,4\}}^{\{4,3\}}\subset\mathscr{P}_{\{3,3,4\}}$, with its cutting plane denoted by $\Sigma$ as the plane $\{z=0\}$. Notice that under the fixed notation described above and the isomorphism $\phi_{\{3,3,4\}}^{\{4,3\}}$ described in \eqref{eq:morphisms}, we have $\mathscr{P}_{\{4,3\}}=\mathscr{S}_{\{3,3,4\}}^{\{4,3\}}\cap\Sigma$. The position of the spheres in $\mathscr{S}_{\{3,3,4\}}^{\{4,3\}}$ relative to $\Sigma$ induces a 2-coloring of $\mathscr{P}_{\{4,3\}}$ by coloring a circle black if it corresponds to a sphere centered above $\Sigma$, and white otherwise (see Figure \ref{fig:BWs}). This coloring is referred to as the \textit{$z$-coloring}, and we shall extended it to the vertices of the tangency graph of $\mathscr{P}_{\{4,3\}}$.

\subsection{The orthocubic shifts}\label{subsec:orthoshifts}
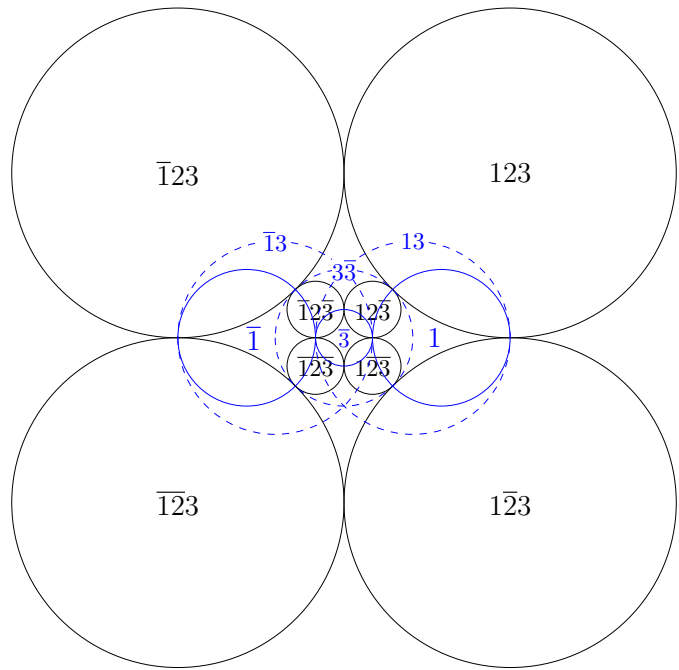
\begin{figure}[H]
			\centering
				\begin{tikzpicture}[scale=1] 

%

%

\node at (.414,.414) [black,fill=white,inner sep=.0pt] { $12\overline{3}$};
\node at (-.414,.414) [black,fill=white,inner sep=.0pt] { $\overline12\overline{3}$};
\node at (.414,-.414) [black,fill=white,inner sep=.0pt] { $1\overline{23}$};
\node at (-.414,-.414) [black,fill=white,inner sep=.0pt] { $\overline{123}$};			
		
		\draw(2.414,2.414) circle (2.414cm) node {\Large $123$};
		\draw (-2.414,2.414) circle (2.414cm) node {\Large$\overline{1}23$};
		\draw (2.414,-2.414) circle (2.414cm) node {\Large$1\overline{2}3$};
		\draw(-2.414,-2.414) circle (2.414cm) node {\Large$\overline{12}3$};
		\draw(.414,.414) circle (.414cm);
		\draw(-.414,.414) circle (.414cm);
		\draw(.414,-.414) circle (.414cm);
		\draw(-.414,-.414) circle (.414cm);

		\draw[blue,dashed] (1.000,0) circle (-1.414cm) node at (1,1.415)  [fill=white, inner sep=0pt] {$13$};
		\draw[blue,dashed] (-1.000,0) circle (-1.414cm) node at (-1,1.415)  [fill=white, inner sep=1pt] {$\overline13$};
		\draw[blue, dashed] (0,0) circle (-1.000cm) node at (0,1)[fill=white, inner sep=0]  { $3\overline3$};

		\draw[blue] (1.414,0) circle (1.00cm) node[xshift=-0.1cm] {\Large$1$};
		\draw[blue,font=\small] (0,0) circle (0.414cm) node {$\overline3$};
		\draw[blue] (-1.414,0) circle (1.00cm) node[xshift=0.1cm] {\Large$\overline1$};

	\end{tikzpicture}
			\caption{ $\mathcal{C}_{\{4,3\}}$ with the walls of the generators of the cubic shifts.}
			\label{fig:dualmirrors}
		\end{figure}

\begin{figure}[H]
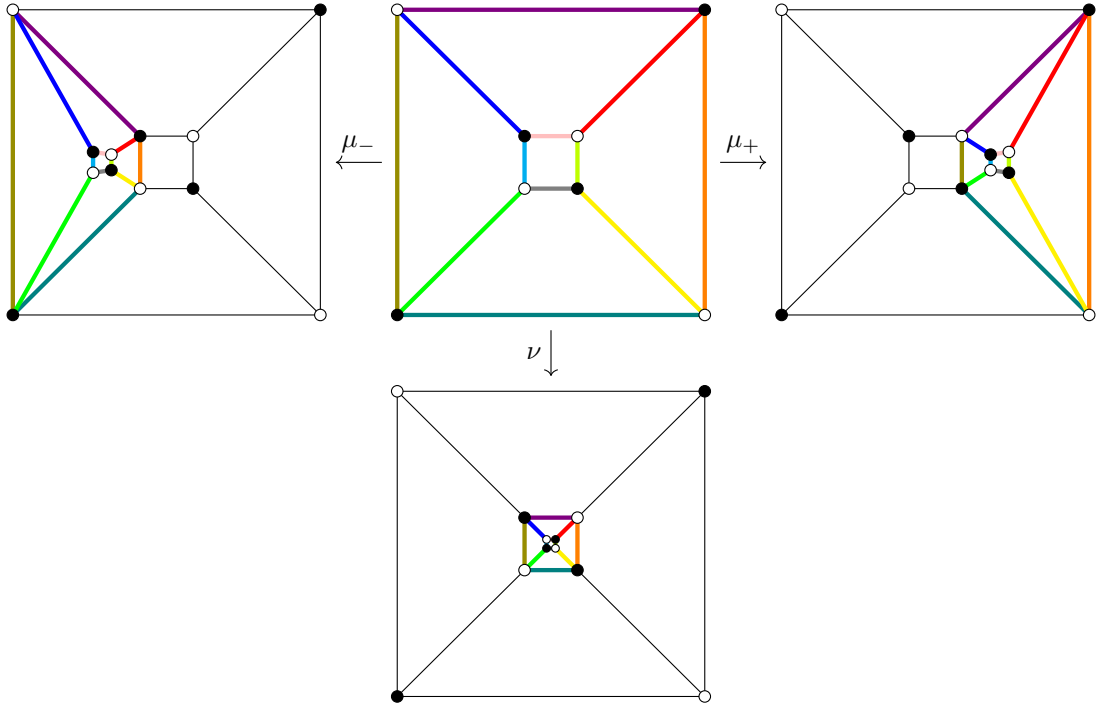

			\centering
			\includestandalone[width=.9\textwidth]{tikzs/orthopushes3}
			\caption{The action of the cubic shifts on the tangency graph of $\mathscr P_{\{4,3\}}$.}
			\label{fig:shifts}
		\end{figure}

We define the \textit{cubic shifts}  as the following parabolic elements:
\begin{align}
	\mu_+=s_1 r_{13} ,&&	\mu_-=s_{\overline 1} r_{\overline13}&&\text{ and }&&\nu:=s_3 r_{3\overline3},
\end{align}
where $r_{13},r_{\overline13},r_{3\overline3},s_1,s_{\overline1},s_{\overline 3}$ are the symmetries of $\Gamma_{\{4,3\}}$ depicted in Figure \ref{fig:dualmirrors}). In Figure \ref{fig:shifts}, we illustrate the action of the cubic shifts on the tangency graph of $\mathscr {P}_{\{4,3\}}$ with the $z$-coloring. Notice that $\mu_{+}$ and $\mu_{-}$ preserves the $z$-coloring while $\nu$ reverses it. We define the \textit{orthocubic shifts} as the symmetries $\widehat{\mu_+}$, $\widehat{\mu_-}$ and $\widehat{\nu}$ of $\mathscr{S}_{\{3,3,4\}}^{\{4,3\}}$ obtained by conjugating the cubic shifts with the morphism induced by \eqref{eq:morphisms}.

\subsection{Orthocubic representations}
	We define an \textit{orthocubic path} $\gamma$ as a polygonal curve in the tangency graph of $\mathscr{S}_{\{3,3,4\}}^{\{4,3\}}$. A \textit{cubic diagram} of $\gamma$ will be its orthogonal projection on $\Sigma$. The orthogonal projection of the tangency graph of $\mathscr{S}_{\{3,3,4\}}^{\{4,3\}}$ on $\Sigma$, is the tangency graph of $\mathscr{P}_{\{4,3\}}$ plus the diagonal edges of each square-face, which join two vertices of same color under the $z$-coloring. Notice that each diagonal edge corresponds to a tangency point of $\mathscr{S}_{\{3,3,4\}}^{\{4,3\}}$ which is above or below $\Sigma$, while each nondiagonal edge corresponds to a tangency point on $\Sigma$. The crossings of any cubic diagram are obtained by the intersection of the two diagonal edges of a same square-face. With the information given by the $z$-coloring, the over/under crossing information can be deduced from the color of the vertices of the diagonal edges  (black=over/white=under).

	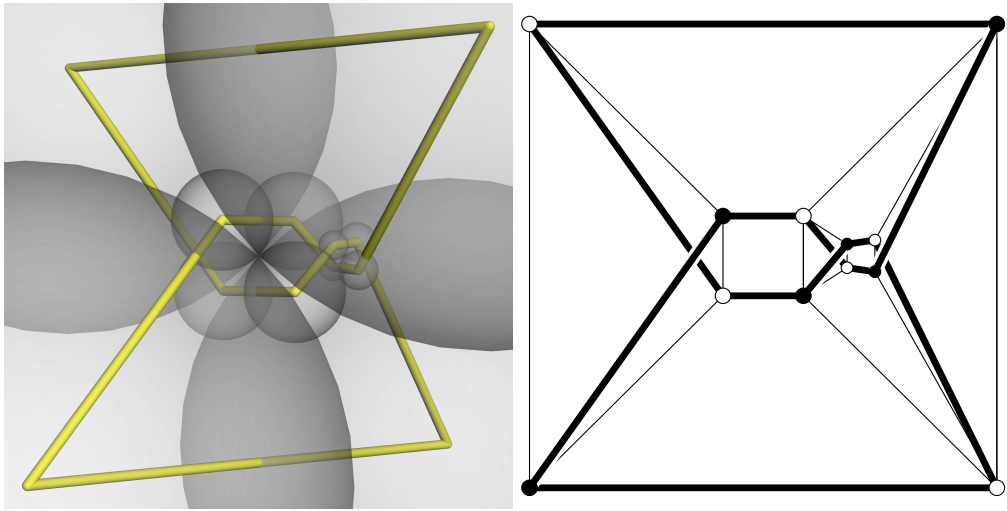
\begin{figure}[H]
		\centering
		\begin{tabular}{cc}
			\includegraphics[width=.42\textwidth,align=c]{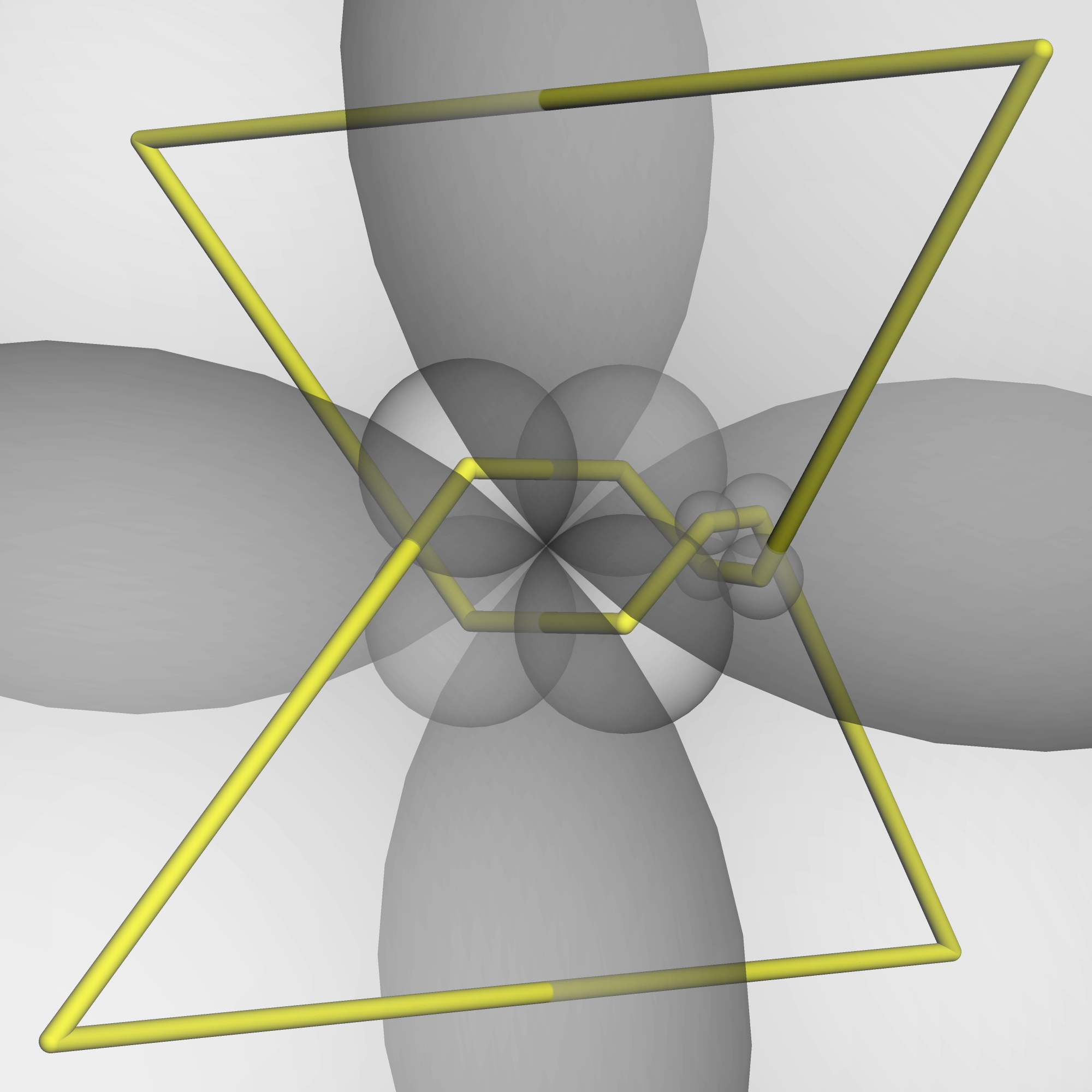}
\begin{tikzpicture}[scale=3]
	\draw \NWin -- \NW;
	\draw \SWin -- \SW;
	\draw \NE -- \SE;
	\draw \NEin -- \SEin;
	\draw \NEin -- \NE;
	\draw \SEin -- \SE;
	\draw \NWr -- \NEr;
	\draw \SWr -- \SEr;
	\draw \NWr -- \SWr;
	\draw \NEr -- \SEr;
	\draw \NE -- \NEr;
	\draw \SE -- \SEr;
	\draw \NEin -- \NWr;
	\draw \SEin -- \SWr;
	\draw[line width=2.5pt]\NE -- \NW;
	\draw[line width=2.5pt]\NWin--\NEin -- \SWr -- \SEr ; 
	\draw[line width=2.5pt] \NW -- \SWin ; 
	\draw[line width=5pt,white] \SEin -- \NWr;
	\draw[line width=2.5pt] \SWin -- \SEin -- \NWr;
	\draw[line width=2.5pt]  \NWr -- \NEr -- \SE -- \SW;
	\draw[line width=5pt,white] \SW -- \NWin;
	\draw[line width=2.5pt] \SW -- \NWin;
	\draw[line width=5pt,white] \SEr -- \NE;
	\draw[line width=2.5pt] \SEr -- \NE; 
	\draw \SW rectangle \NE;
	\draw \SWin rectangle \NEin;
	\node  at \NE [circle,draw=black,fill=black,inner sep=0pt,minimum size=.2cm]  {};
	\node  at \NW [circle,draw=black,fill=white,inner sep=0pt,minimum size=.2cm]  {};
	\node  at \SW [circle,draw=black,fill=black,inner sep=0pt,minimum size=.2cm]  {};
	\node  at \SE [circle,draw=black,fill=white,inner sep=0pt,minimum size=.2cm]  {};
	\node  at \NEin [circle,draw=black,fill=white,inner sep=0pt,minimum size=.2cm]  {};
	\node  at \NWin [circle,draw=black,fill=black,inner sep=0pt,minimum size=.2cm]  {};
	\node  at \SWin [circle,draw=black,fill=white,inner sep=0pt,minimum size=.2cm]  {};
	\node  at \SEin [circle,draw=black,fill=black,inner sep=0pt,minimum size=.2cm]  {};
	\node  at \NWr [circle,draw=black,fill=black,inner sep=0pt,minimum size=.15cm]  {};
	\node  at \NEr [circle,draw=black,fill=white,inner sep=0pt,minimum size=.15cm]  {};
	\node  at \SEr [circle,draw=black,fill=black,inner sep=0pt,minimum size=.15cm]  {};
	\node at \SWr  [circle,draw=black,fill=white,inner sep=0pt,minimum size=.15cm]  {};
\end{tikzpicture} &
		\end{tabular}
		\caption{(Left) An orthocubic representation of the trefoil knot and its corresponding cubic diagram (right).}\label{fig:orthotrefoil}
	\end{figure} 
 We define an \textit{orthocubic representation} of a link $L$ as a collection of disjoint closed orthocubic paths isotopically equivalent to $L$. Every orthocubic representation induces a necklace representation in $\mathscr{S}_{\{3,3,4\}}^{\{4,3\}}$. In Figure \ref{fig:orthotrefoil}, we show an orthocubic representation of the trefoil knot, and its corresponding cubic diagram.	Let $T$ be the tetrahedron whose vertices are the centers of the spheres $\{S_{\overline1},S_{\overline2},S_3,S_4\}\subset\mathcal{S}_{\{3,3,4\}}$. We define an \textit{orthocubic tangle} as a tangle $(T,\boxed{t})$ where $\boxed{t}$ is a collection $\{\gamma_1,\gamma_2,\ldots,\gamma_m\}$ of $m\ge2$ disjoint orthocubic paths contained in $T$, satisfying that the endpoints of $\gamma_1$ and $\gamma_2$ lie in the corners of $T$, and the rest of the orthocubic paths are closed. We now have all the necessary framework to prove the main result of this section.
			 
	\begin{thm}\label{thm:orthocubicrep} Every algebraic link admits a necklace representation contained in $\mathscr{S}_{\{3,3,4\}}^{\{4,3\}}$.
	\end{thm}
	\begin{proof}
		We begin by constructing the orthocubic analogs of the building blocks $(\T_0)-(\T_5)$ of algebraic tangles defined in Section \ref{sec:preliminaries}.
				
		\begin{enumerate}[label=($\O_\arabic*)$]
			\setcounter{enumi}{-1}
			\item The \textit{orthocubic elementary tangles} are described in Figure \ref{fig:orthoelem}.
			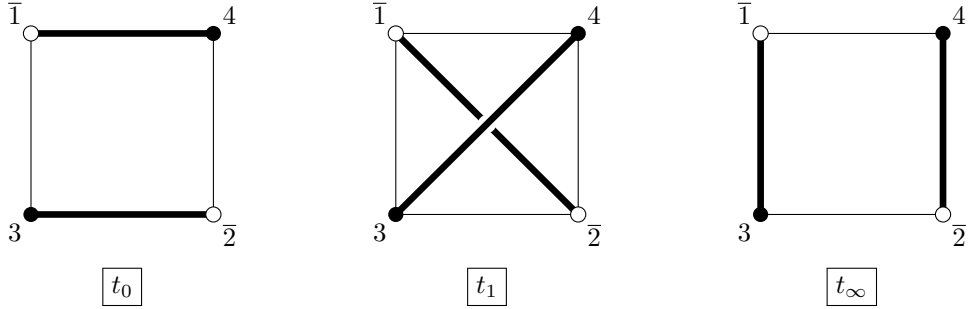
\begin{figure}[H]
			\centering
	\begin{tikzpicture}[scale=1.2]
		\begin{scope}[xshift=-4cm]
			\draw \SW rectangle \NE;
			
			\draw[line width=2.5pt] \NW -- \NE;
			\draw[line width=2.5pt] \SW -- \SE;
			
			\draw[fill=black] \NE circle (0cm) node[anchor= south west] {$4$};
			\draw[draw=black,fill=white] \NW circle (0cm) node[anchor= south east] {$\overline1$};
			\draw[fill=black] \SW circle (0cm) node[anchor=north  east] {$3$};
			\draw[draw=black,fill=white] \SE circle (0cm) node[anchor=north west] {$\overline2$};
			
			\node at \NE [circle,fill=black,inner sep=0pt,minimum size=.2cm]  {};
			\node at \NW [circle,draw=black,fill=white,inner sep=0pt,minimum size=.2cm]  {};
			\node at \SW [circle,fill=black,inner sep=0pt,minimum size=.2cm]  {};
			\node at \SE [circle,draw=black,fill=white,inner sep=0pt,minimum size=.2cm]  {};
			
			\node at (0,-1.8) {$\boxed{t_0}$};
		\end{scope}
		
		
		\draw \SW rectangle \NE;
		

		\draw[line width=2.5pt] \NW -- \SE;
		\fill[white] circle (.1cm);
		\draw[line width=2.5pt] \SW -- \NE;
		
		\draw[fill=black] \NE circle (0cm) node[anchor= south west] {$4$};
		\draw[draw=black,fill=white] \NW circle (0cm) node[anchor= south east] {$\overline1$};
		\draw[fill=black] \SW circle (0cm) node[anchor=north  east] {$3$};
		\draw[draw=black,fill=white] \SE circle (0cm) node[anchor=north west] {$\overline2$};
		
		\node at \NE [circle,fill=black,inner sep=0pt,minimum size=.2cm]  {};
		\node at \NW [circle,draw=black,fill=white,inner sep=0pt,minimum size=.2cm]  {};
		\node at \SW [circle,fill=black,inner sep=0pt,minimum size=.2cm]  {};
		\node at \SE [circle,draw=black,fill=white,inner sep=0pt,minimum size=.2cm]  {};
		
		\node at (0,-1.8) {$\boxed{t_1}$};
		
		\begin{scope}[xshift=4cm]
			\draw \SW rectangle \NE;
			
			\draw[line width=2.5pt] \SW -- \NW;
			\draw[line width=2.5pt] \NE -- \SE;
			
			\draw[fill=black] \NE circle (.0cm) node[anchor= south west] {$4$};
			\draw[draw=black,fill=white] \NW circle (.0cm) node[anchor= south east] {$\overline1$};
			\draw[fill=black] \SW circle (.0cm) node[anchor=north  east] {$3$};
			\draw[draw=black,fill=white] \SE circle (.0cm) node[anchor=north west] {$\overline2$};
			
			\node at \NE [circle,fill=black,inner sep=0pt,minimum size=.2cm]  {};
			\node at \NW [circle,draw=black,fill=white,inner sep=0pt,minimum size=.2cm]  {};
			\node at \SW [circle,fill=black,inner sep=0pt,minimum size=.2cm]  {};
			\node at \SE [circle,draw=black,fill=white,inner sep=0pt,minimum size=.2cm]  {};
			
			\node at (0,-1.8) {$\boxed{t_\infty}$};
		\end{scope}
		
	\end{tikzpicture}   
			\caption{The orthocubic elementary tangles.}\label{fig:orthoelem}
			\end{figure} 				
			\item The \textit{orthocubic flip} $ F_\O\boxed{t}:= \widehat{r}_{12}\boxed{t}$, where $\widehat{r}_{12}\in \Gamma_{\{3,3,4\}}$ is the reflection on the plane $\{y=x\}$.
			\item The \textit{orthocubic mirror} $-\boxed{t}:=\widehat \nu\boxed{t}$ plus the edges $\{( 1,\overline2),(\overline1,2),( 3,\overline4),(\overline3,4)\}$.
			\begin{figure}[H]
	\begin{tikzpicture}[scale=1.2]
		\draw[fill=mylgray]  \SW rectangle \NE;
		\node at \NE [circle,fill=black,inner sep=0pt,minimum size=.2cm]  {};
		\node at \NW [circle,draw=black,fill=white,inner sep=0pt,minimum size=.2cm]  {};
		\node at \SW [circle,fill=black,inner sep=0pt,minimum size=.2cm]  {};
		\node at \SE [circle,draw=black,fill=white,inner sep=0pt,minimum size=.2cm]  {};
		
		\node at \NE [anchor=south west]  {$4$};
		\node at \NW [anchor=south east]  {$\overline1$};
		\node at \SE [anchor=north west]  {$\overline2$};
		\node at \SW [anchor=north east]  {$3$};
		\node[] {{\fontsize{40}{60} $t$}};
		
		\node at (0,-1.8) {$\boxed t$};
		
%
		
		\begin{scope}[xshift=4cm]
			
			\draw[fill=mylgray]  \SW rectangle \NE;
			
			\node at \NE [circle,fill=black,inner sep=0pt,minimum size=.2cm]  {};
			\node at \NW [circle,draw=black,fill=white,inner sep=0pt,minimum size=.2cm]  {};
			\node at \SW [circle,fill=black,inner sep=0pt,minimum size=.2cm]  {};
			\node at \SE [circle,draw=black,fill=white,inner sep=0pt,minimum size=.2cm]  {};
			
			\node at \NE [anchor=south west]  {$4$};
			\node at \NW [anchor=south east]  {$\overline1$};
			\node at \SE [anchor=north west]  {$\overline2$};
			\node at \SW [anchor=north east]  {$3$};
			
			\node[] {{\fontsize{40}{60} \rotatebox{-90}{\reflectbox{$t$}}}};
			\draw[dashed] \SW -- \NE;
			
			\node at (0,-1.8) {$F_\O\boxed t$};
		\end{scope}
	
		---(-t)
	\begin{scope}[xshift=8cm]
		\draw \SW rectangle \NE;
		
		\draw[fill=mylgray] \SWin rectangle \NEin;
		
		\draw[line width=2.5pt] \NEin -- \NE;
		\draw[line width=2.5pt] \NWin --\NW;
		\draw[line width=2.5pt] \SEin -- \SE;
		\draw[line width=2.5pt] \SWin --\SW;
		
		\node at \NE [circle,fill=black,inner sep=0pt,minimum size=.2cm]  {};
		\node at \NW [circle,draw=black,fill=white,inner sep=0pt,minimum size=.2cm]  {};
		\node at \SW [circle,fill=black,inner sep=0pt,minimum size=.2cm]  {};
		\node at \SE [circle,draw=black,fill=white,inner sep=0pt,minimum size=.2cm]  {};
		
		\node at \NEin  [circle,draw=black,fill=white,inner sep=0pt,minimum size=.15cm]  {};
		\node at \NWin [circle,fill=black,inner sep=0pt,minimum size=.15cm]  {};
		\node at \SWin  [circle,draw=black,fill=white,inner sep=0pt,minimum size=.15cm]  {};
		\node at \SEin [circle,fill=black,inner sep=0pt,minimum size=.15cm]  {};
		
						\node at \NE [anchor=south west]  {$4$};
\node at \NW [anchor=south east]  {$\overline1$};
\node at \SE [anchor=north west]  {$\overline2$};
\node at \SW [anchor=north east]  {$3$};
\node at \NEin [anchor=south,inner sep=5pt ]  {\scriptsize $\overline3$};
\node at \NWin [anchor=south,inner sep=5pt  ]  {\scriptsize $2$};
\node at \SEin [anchor=north,inner sep=6pt  ]  {\scriptsize $1$};
\node at \SWin [anchor=north,inner sep=4pt  ]  {\scriptsize $\overline4$};
		
		\node[] {{$t$}};
		
		\node at (0,-1.8) {$-\boxed t$};
		
	\end{scope}

	\end{tikzpicture}  
	
					\caption{The orthocubic flip and the orthocubic mirror.}
			\end{figure}
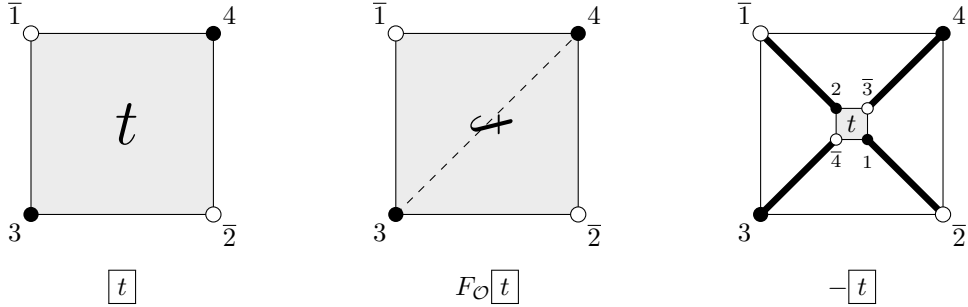 				
	
			\item The \textit{orthocubic addition} $\boxed{t'}+ \boxed{t}:=\widehat\mu_{-}\boxed{t'}\cup \widehat\mu_+\boxed{t}$ plus the edges $\{ (1,\overline4),(2,\overline3)\}$.
			\item The \textit{orthocubic half-twists} $H_\O^+\boxed{t}:=\boxed{t_1}+\boxed{t}$ and  $H^-_{\O}\boxed{t}:=-\boxed{t_{1}}+ \boxed{t}$.
			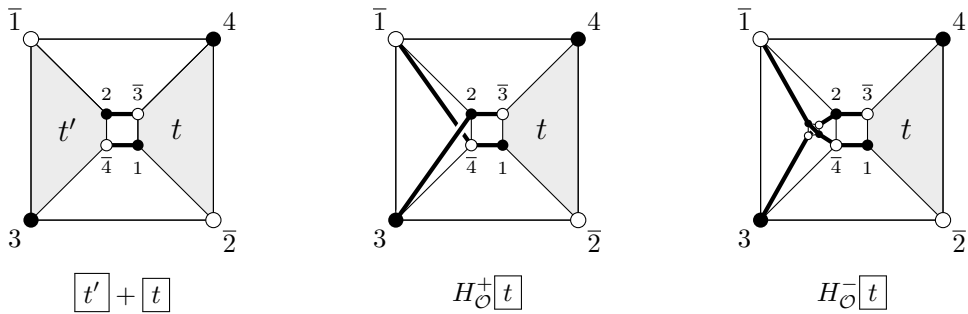
\begin{figure}[H]
				\centering
						\renewcommand{\thi}{.7mm}
	\begin{tikzpicture}[scale=1.2]

\begin{scope}
	
	\draw[ultra thick] \NEin -- \NWin;
	\draw[ultra thick] \SWin -- \SEin;
	
	\draw[fill=mylgray] \SWin -- \NWin --\NW -- \SW;
	\draw[fill=mylgray] \SEin -- \NEin  --\NE -- \SE;
	\draw \SW -- \SE;
	\draw \NW -- \NE;
	
	\draw \NEin  -- \NE;
	\draw \NWin -- \NW;
	\draw \SEin -- \SE;
	\draw \SWin -- \SW;
	
	\draw \SW rectangle \NE;
	
	\draw \SWin rectangle \NEin ;
	
	\node at \NE [circle,fill=black,inner sep=0pt,minimum size=.2cm]  {};
	\node at \NW [circle,draw=black,fill=white,inner sep=0pt,minimum size=.2cm]  {};
	\node at \SW [circle,fill=black,inner sep=0pt,minimum size=.2cm]  {};
	\node at \SE [circle,draw=black,fill=white,inner sep=0pt,minimum size=.2cm]  {};
	
	\node (NEin) at \NEin  [circle,draw=black,fill=white,inner sep=0pt,minimum size=.15cm]  {};
	\node (NWin) at \NWin [circle,fill=black,inner sep=0pt,minimum size=.15cm]  {};
	\node (SWin) at \SWin  [circle,draw=black,fill=white,inner sep=0pt,minimum size=.15cm]  {};
	\node (SEin) at \SEin [circle,fill=black,inner sep=0pt,minimum size=.15cm]  {};
	
	\node at (-.6,0) {{\Large $t'$}};
	\node at (.6,0) {{\Large $t$}};
	
				\node at \NE [anchor=south west]  {$4$};
				\node at \NW [anchor=south east]  {$\overline1$};
				\node at \SE [anchor=north west]  {$\overline2$};
				\node at \SW [anchor=north east]  {$3$};
\node at \NEin [anchor=south,inner sep=5pt ]  {\scriptsize $\overline3$};
\node at \NWin [anchor=south,inner sep=5pt  ]  {\scriptsize $2$};
\node at \SEin [anchor=north,inner sep=6pt  ]  {\scriptsize $1$};
\node at \SWin [anchor=north,inner sep=4pt  ]  {\scriptsize $\overline4$};
				
	\node at (0,-1.8) {$\boxed {t'} + \boxed {t}$};
	
\end{scope}
		
				\begin{scope}[xshift=4cm]
						\filldraw[mylgray] \SEin --\NEin --\NE -- \SE;
						\draw \SW rectangle \NE;
						

						\draw[ultra thick] \NW-- \SWin -- \SEin;
						\draw[white, line width=4pt] \SW -- \NWin --\NEin;
						\draw[ultra thick] \SW -- \NWin --\NEin;
						
						\draw \NEin -- \NE;
						\draw \NWin -- \NW;
						\draw \SEin -- \SE;
						\draw \SWin -- \SW;
						
						\draw \SW rectangle \NE;
						
						\draw \SWin rectangle \NEin;
						
						\node  at \NE [circle,fill=black,inner sep=0pt,minimum size=.2cm]  {};
						\node  at \NW [circle,draw=black,fill=white,inner sep=0pt,minimum size=.2cm]  {};
						\node  at \SW [circle,fill=black,inner sep=0pt,minimum size=.2cm]  {};
						\node  at \SE [circle,draw=black,fill=white,inner sep=0pt,minimum size=.2cm]  {};
						
						\node at \NEin  [circle,draw=black,fill=white,inner sep=0pt,minimum size=.15cm]  {};
						\node  at \NWin [circle,fill=black,inner sep=0pt,minimum size=.15cm]  {};
						\node  at \SWin  [circle,draw=black,fill=white,inner sep=0pt,minimum size=.15cm]  {};
						\node  at \SEin [circle,fill=black,inner sep=0pt,minimum size=.15cm]  {};
						
						\node at \NE [anchor=south west]  {$4$};
						\node at \NW [anchor=south east]  {$\overline1$};
						\node at \SE [anchor=north west]  {$\overline2$};
						\node at \SW [anchor=north east]  {$3$};
\node at \NEin [anchor=south,inner sep=5pt ]  {\scriptsize $\overline3$};
\node at \NWin [anchor=south,inner sep=5pt  ]  {\scriptsize $2$};
\node at \SEin [anchor=north,inner sep=6pt  ]  {\scriptsize $1$};
\node at \SWin [anchor=north,inner sep=4pt  ]  {\scriptsize $\overline4$};
						
						\node at (.6,0) {{\Large $t$}};
						\node at (0,-1.8) {$H_\O^+\boxed t$};
					\end{scope}

				\begin{scope}[xshift=8cm]
						\filldraw[mylgray] \SEin --\NEin --\NE -- \SE;
						\draw \SW rectangle \NE;
						
						
						\draw[black,ultra thick] \SW --\SWl--\NEl-- \NWin--\NEin;
						\draw[white,line width=2.5pt] \NW --\NWl--\SEl-- \SWin--\SEin;
						\draw[black,ultra thick] \NW --\NWl--\SEl-- \SWin--\SEin;
						
						\draw \NEin -- \NE;
						\draw \NWin -- \NW;
						\draw \SEin -- \SE;
						\draw \SWin -- \SW;
						
						\draw \NWl--\NEl--\SEl--\SWl--\NWl;
						
						\draw \SW rectangle \NE;
						
						\draw \SWin rectangle \NEin;
						
						\node  at \NE [circle,fill=black,inner sep=0pt,minimum size=.2cm]  {};
						\node  at \NW [circle,draw=black,fill=white,inner sep=0pt,minimum size=.2cm]  {};
						\node  at \SW [circle,fill=black,inner sep=0pt,minimum size=.2cm]  {};
						\node  at \SE [circle,draw=black,fill=white,inner sep=0pt,minimum size=.2cm]  {};
						
						\node at \NEin  [circle,draw=black,fill=white,inner sep=0pt,minimum size=.15cm]  {};
						\node at \NWin [circle,fill=black,inner sep=0pt,minimum size=.15cm]  {};
						\node at \SWin  [circle,draw=black,fill=white,inner sep=0pt,minimum size=.15cm]  {};
						\node at \SEin [circle,fill=black,inner sep=0pt,minimum size=.15cm]  {};
						
						\node at \NEl  [circle,draw=black,fill=white,inner sep=0pt,minimum size=.1cm]  {};
						\node at \NWl [circle,fill=black, inner sep=0pt,minimum size=.1cm]  {};
						\node at \SWl [circle,draw=black,fill=white,inner sep=0pt,minimum size=.1cm]  {};
						\node at \SEl [circle,fill=black, inner sep=0pt,minimum size=.1cm]  {};
						
						\node at \NE [anchor=south west]  {$4$};
						\node at \NW [anchor=south east]  {$\overline1$};
						\node at \SE [anchor=north west]  {$\overline2$};
						\node at \SW [anchor=north east]  {$3$};
\node at \NEin [anchor=south,inner sep=5pt ]  {\scriptsize $\overline3$};
\node at \NWin [anchor=south,inner sep=5pt  ]  {\scriptsize $2$};
\node at \SEin [anchor=north,inner sep=6pt  ]  {\scriptsize $1$};
\node at \SWin [anchor=north,inner sep=4pt  ]  {\scriptsize $\overline4$};

						\node at (.6,0) {{\Large $t$}};
						\node at (0,-1.8) {$H_\O^-\boxed t$};
					\end{scope}
	\end{tikzpicture}  
	
				\caption{The orthocubic addition and the orthocubic half-twists.} 
			\end{figure}			
			\item 
				\renewcommand{\thi}{.7mm}
					\newcommand{\lwidth}{.5mm}
					\newcommand{\edgeinb}{.4}
					\newcommand{\NEinc}{(\edgeinb,\edgeinb)}
					\newcommand{\NWinc}{(-\edgeinb,\edgeinb)}
					\newcommand{\SWinc}{(-\edgeinb,-\edgeinb)}
					\newcommand{\SEinc}{(\edgeinb,-\edgeinb)}
					
			The  \textit{orthocubic tangle closures} are described in Figure \ref{fig:tangclose}.
				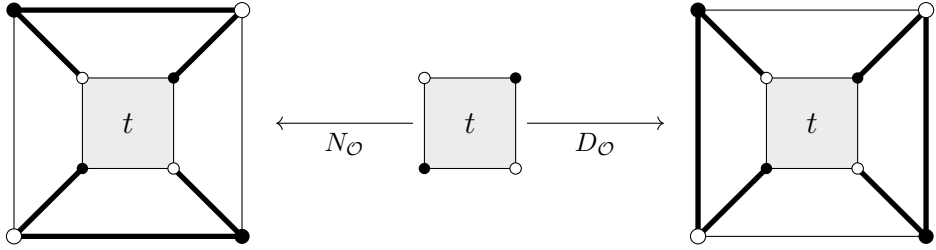
\begin{figure}[H]
					\centering
	\begin{tikzpicture}[scale=1.5]
		\draw[->] (0.5,0) -> (1.7,0) node [below,pos=.5] {$D_{\mathcal O}$};
		\draw[->] (-0.5,0) -> (-1.7,0) node [below,pos=.5] {$N_{\mathcal O}$};		
		\begin{scope}[xshift=0cm]

			\draw[fill=mylgray] \SWinc rectangle \NEinc;

			\node  at \NEinc  [circle,fill=black,inner sep=0pt,minimum size=.15cm]  {};
			\node  at \NWinc [circle,draw=black,fill=white,inner sep=0pt,minimum size=.15cm]  {};
			\node  at \SWinc  [circle,fill=black,inner sep=0pt,minimum size=.15cm]  {};
			\node  at \SEinc [circle,draw=black,fill=white,inner sep=0pt,minimum size=.15cm]  {};
			
			\node[] {{\Large $t$}};
			
			
			%
		\end{scope}

		\begin{scope}[xshift=-3cm]

			\draw \SW rectangle \NE;
			\draw[fill=mylgray] \SWinc rectangle \NEinc;
			
			\draw[line width=\thi] \NEinc -- \NE -- \NW -- \NWinc;
			\draw[line width=\thi] \SEinc -- \SE -- \SW-- \SWinc;
			
			\node  at \NE [circle,draw=black,fill=white,inner sep=0pt,minimum size=.2cm]  {};
			\node at \NW [circle,fill=black,inner sep=0pt,minimum size=.2cm]  {};
			\node  at \SW [circle,draw=black,fill=white,inner sep=0pt,minimum size=.2cm]  {};
			\node  at \SE [circle,fill=black,inner sep=0pt,minimum size=.2cm]  {};
			
			\node  at \NEinc  [circle,fill=black,inner sep=0pt,minimum size=.15cm]  {};
			\node  at \NWinc [circle,draw=black,fill=white,inner sep=0pt,minimum size=.15cm]  {};
			\node  at \SWinc  [circle,fill=black,inner sep=0pt,minimum size=.15cm]  {};
			\node  at \SEinc [circle,draw=black,fill=white,inner sep=0pt,minimum size=.15cm]  {};
			
			\node[] {{\Large $t$}};
			%

		\end{scope}
		
			

		\begin{scope}[xshift=3cm]
			
			\draw \SW rectangle \NE;
			\draw[fill=mylgray] \SWinc rectangle \NEinc;
			
			\draw[line width=\thi] \NEinc -- \NE -- \SE -- \SEinc;
			\draw[line width=\thi] \SWinc -- \SW -- \NW-- \NWinc;
			
			\node  at \NE [circle,draw=black,fill=white,inner sep=0pt,minimum size=.2cm]  {};
			\node at \NW [circle,fill=black,inner sep=0pt,minimum size=.2cm]  {};
			\node  at \SW [circle,draw=black,fill=white,inner sep=0pt,minimum size=.2cm]  {};
			\node  at \SE [circle,fill=black,inner sep=0pt,minimum size=.2cm]  {};
			
			\node  at \NEinc  [circle,fill=black,inner sep=0pt,minimum size=.15cm]  {};
			\node  at \NWinc [circle,draw=black,fill=white,inner sep=0pt,minimum size=.15cm]  {};
			\node  at \SWinc  [circle,fill=black,inner sep=0pt,minimum size=.15cm]  {};
			\node  at \SEinc [circle,draw=black,fill=white,inner sep=0pt,minimum size=.15cm]  {};
			
			\node[] {{\Large $t$}};

		\end{scope}
		
	\end{tikzpicture}  
					\caption{The orthocubic tangle closures.}
					\label{fig:tangclose}
				\end{figure} 
			
		\end{enumerate}
				
		For every $i=0,\ldots,5$, the orthocubic building block $(\O_i)$ is  isotopically equivalent to the building block $(\T_i)$ of algebraic tangles. Therefore, we can mimic the Conway's algorithm to define an \textit{orthocubic} rational tangle isotopically equivalent to $t[a_1,\cdots,a_n]$ by
			\begin{align}\label{eq!orthoConway}
				t_\O[a_1,\cdots,a_n]:=H_\O^{a_1}F_\O\cdots H_\O^{a_n}F_\O\boxed{t_\infty}
			\end{align}
		By combining orthocubic addition, flip and closure of orthocubic rational tangles we obtain orthocubic representations of every algebraic link.
	\end{proof}

\subsection{Improvement of the upper bound of the ball number}				
	The orthocubic Conway's algorithm can be slightly adapted in order to attain the upper bound of Theorem \ref{thm:rationalball}. For every $a_1\ge0$, $a_2,\ldots,a_n>0$, we define the \textit{reduced} orthocubic Conway's algorithm $\widetilde t_\O[a_1,\cdots,a_n]$ by 
	\begin{align}
		\widetilde t_\O[a_1,\cdots,a_n]:=H_\O^{a_1}F_\O\cdots H_\O^{a_n-1}\boxed{t_1}
	\end{align}
	Clearly, for every $a_1\ge0$, $a_2,\ldots,a_n>0$, we have $t_\O[a_1,\cdots,a_n]\simeq \widetilde t_\O[a_1,\cdots,a_n]$.	
	\begin{thm}\label{thm:rationalball} Let $L$ be an algebraic link obtained by the closure of the algebraic tangle  $t_{p_1/q_1}+\cdots+t_{p_m/q_m}$ where all the $p_i/q_i$ have same sign. Then,
		$\mathrm {ball}(L)\leq 4\mathrm{cr}(L).$
	\end{thm}

	\begin{proof} 
	Let $L$ be an algebraic link made by the closure $N(t)$ where $t=t_{p_1/q_1}+\cdots+t_{p_m/q_m}.$
	The condition that all $p_i/q_i$ have the same sign implies that we have alternating diagram of $L$ induced by the closure of $t$. Thus, the crossing number of $L$ is equal to the sum of the crossing numbers of each $t_{p_i/q_i}$.  Without loss of generality, we can consider that all $p_i/q_i$ are positive.
	\medskip
	
	 For every $p_i/q_i$ with positive continued fraction $[a_1,\cdots,a_n]$, let $\boxed{t_{p_i/q_i}}:=\widetilde t_\O[a_1,\cdots,a_n]$. Since $F_\O$ does not change the number of the spheres of the corresponding necklace, and $H_\O^+$ increases it by $4$, we have that the number of spheres in $\boxed{t_{p_i/q_i}}$ is equal to $4(a_1+\ldots+a_{n})=4\mathrm{cr}(t_{p_i/q_i}).$ 
	 \medskip
	 
	 Now, let $\boxed t$ be the orthocubic tangle made by the orthocubic additions $\boxed{t_{p_1/q_1}}+\cdots+\boxed{t_{p_m/q_m}}$. By the equivalence between the orthocubic and tangle operations we have that $\boxed t\simeq t$. Since the function counting the number of spheres is additive for the orthocubic addition, we have that the number of spheres of $\boxed t$ is equal to $4\mathrm{cr}(t_{p_1/q_1})+\cdots+4\mathrm{cr}(t_{p_m/q_m})=4\mathrm{cr}(L)$.
	 \medskip
	 
	  Finally, since the exterior edges $(\overline1,4)$ and $(\overline2,3)$ are not included in any orthocubic tangle obtained after applying an orthocubic addition, we can use these edges to close $\boxed{t}$, and thus obtain a necklace representation of $L$ with $4\mathrm{cr}(L)$ spheres. 
	\end{proof}

\subsection{No tightness for non-alternating links}				
	The family of algebraic links considered in Theorem \ref{thm:rationalball} contains all the rational links and other well-known families as the \textit{Montesinos links} with positive coefficients. These are the links obtained by the closure of $t_{p_1/q_1}+\cdots+t_{p_n/q_n}+t_r$
	with $p_i/q_i>0$ and $r\ge0$. If $r=0$ and every $p_i=1$, then we obtain the \textit{Pretzel link} $P(q_1,\ldots,q_n)$.
\medskip

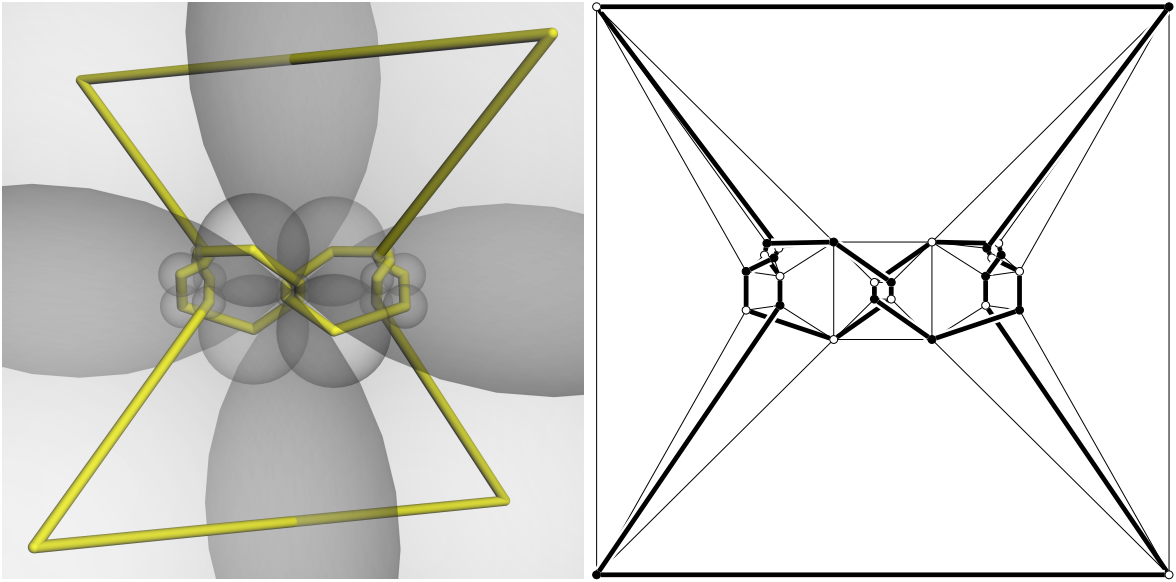
\begin{figure}[H]
		\centering
		\includegraphics[width=.48\textwidth,align=c]{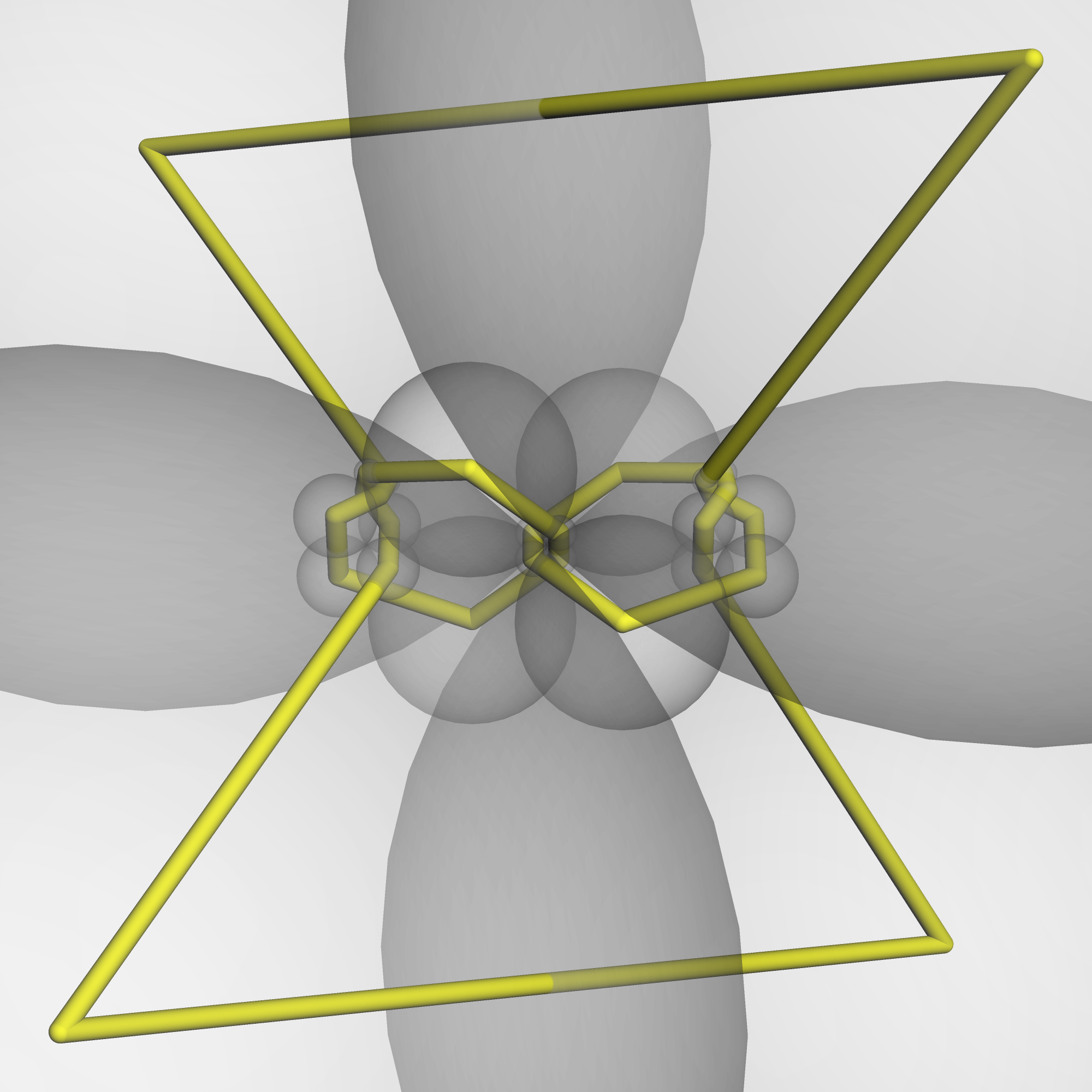}	
\begin{tikzpicture}[scale=3]
\draw[line width=.06mm] (0.4142,0.4142) -- (0.4142,-0.4142);
\draw[line width=.06mm] (0.4142,0.4142) -- (0.8770,0.3611);
\draw[line width=.06mm] (0.4142,0.4142) -- (2.414,2.414);
\draw[line width=.06mm] (0.4142,0.4142) -- (0.07107,0.07107);
\draw[line width=.06mm] (0.4142,0.4142) -- (-0.4142,0.4142);
\draw[line width=.06mm] (0.4142,0.4142) -- (0.8673,0.1239);
\draw[line width=.06mm] (-0.07107,0.07107) -- (-0.07107,-0.07107);
\draw[line width=.06mm] (-0.07107,0.07107) -- (0.07107,0.07107);
\draw[line width=.06mm] (-0.07107,0.07107) -- (-0.4142,0.4142);
\draw[line width=.06mm] (-0.07107,-0.07107) -- (-0.4142,-0.4142);
\draw[line width=.06mm] (-0.07107,-0.07107) -- (0.07107,-0.07107);
\draw[line width=.06mm] (0.4142,-0.4142) -- (-0.4142,-0.4142);
\draw[line width=.06mm] (0.4142,-0.4142) -- (0.07107,-0.07107);
\draw[line width=.06mm] (0.4142,-0.4142) -- (2.414,-2.414);
\draw[line width=.06mm] (0.4142,-0.4142) -- (0.8673,-0.1239);
\draw[line width=.06mm] (1.153,-0.1647) -- (1.153,0.1647);
\draw[line width=.06mm] (1.153,-0.1647) -- (2.414,-2.414);
\draw[line width=.06mm] (1.153,-0.1647) -- (0.8673,-0.1239);
\draw[line width=.06mm] (1.153,0.1647) -- (2.414,2.414);
\draw[line width=.06mm] (1.153,0.1647) -- (0.8673,0.1239);
\draw[line width=.06mm] (1.153,0.1647) -- (0.9982,0.3038);
\draw[line width=.06mm] (0.9185,0.2795) -- (0.8770,0.3611);
\draw[line width=.06mm] (0.9185,0.2795) -- (0.8673,0.1239);
\draw[line width=.06mm] (0.9185,0.2795) -- (0.9982,0.3038);
\draw[line width=.06mm] (0.8770,0.3611) -- (0.9779,0.4027);
\draw[line width=.06mm] (2.414,2.414) -- (-2.414,2.414);
\draw[line width=.06mm] (2.414,2.414) -- (2.414,-2.414);
\draw[line width=.06mm] (2.414,2.414) -- (0.9779,0.4027);
\draw[line width=.06mm] (-2.414,2.414) -- (-1.153,0.1647);
\draw[line width=.06mm] (-2.414,2.414) -- (-0.4142,0.4142);
\draw[line width=.06mm] (-2.414,2.414) -- (-0.9779,0.4027);
\draw[line width=.06mm] (-2.414,2.414) -- (-2.414,-2.414);
\draw[line width=.06mm] (-0.8770,0.3611) -- (-0.9185,0.2795);
\draw[line width=.06mm] (-0.8770,0.3611) -- (-0.4142,0.4142);
\draw[line width=.06mm] (-0.8770,0.3611) -- (-0.9779,0.4027);
\draw[line width=.06mm] (-0.9185,0.2795) -- (-0.9982,0.3038);
\draw[line width=.06mm] (-0.9185,0.2795) -- (-0.8673,0.1239);
\draw[line width=.06mm] (-1.153,0.1647) -- (-1.153,-0.1647);
\draw[line width=.06mm] (-1.153,0.1647) -- (-0.9982,0.3038);
\draw[line width=.06mm] (-1.153,0.1647) -- (-0.8673,0.1239);
\draw[line width=.06mm] (-1.153,-0.1647) -- (-0.8673,-0.1239);
\draw[line width=.06mm] (-1.153,-0.1647) -- (-2.414,-2.414);
\draw[line width=.06mm] (-0.4142,-0.4142) -- (-0.4142,0.4142);
\draw[line width=.06mm] (-0.4142,-0.4142) -- (-0.8673,-0.1239);
\draw[line width=.06mm] (-0.4142,-0.4142) -- (-2.414,-2.414);
\draw[line width=.06mm] (0.07107,-0.07107) -- (0.07107,0.07107);
\draw[line width=.06mm] (-0.4142,0.4142) -- (-0.8673,0.1239);
\draw[line width=.06mm] (-0.9779,0.4027) -- (-0.9982,0.3038);
\draw[line width=.06mm] (-0.8673,0.1239) -- (-0.8673,-0.1239);
\draw[line width=.06mm] (-2.414,-2.414) -- (2.414,-2.414);
\draw[line width=.06mm] (0.8673,-0.1239) -- (0.8673,0.1239);
\draw[line width=.06mm] (0.9982,0.3038) -- (0.9779,0.4027);

\node (1) at (0.4142,0.4142) [circle,draw=black,fill=white,inner sep=0pt,minimum size=\rad] {};
\node (2) at (-0.07107,0.07107) [circle,draw=black,fill=white,inner sep=0pt,minimum size=\rad] {};
\node (3) at (-0.07107,-0.07107) [circle,draw=black,fill=black,inner sep=0pt,minimum size=\rad] {};
\node (4) at (0.4142,-0.4142) [circle,draw=black,fill=black,inner sep=0pt,minimum size=\rad] {};
\node (5) at (1.153,-0.1647) [circle,draw=black,fill=black,inner sep=0pt,minimum size=\rad] {};
\node (6) at (1.153,0.1647) [circle,draw=black,fill=white,inner sep=0pt,minimum size=\rad] {};
\node (7) at (0.9185,0.2795) [circle,draw=black,fill=white,inner sep=0pt,minimum size=\rad] {};
\node (8) at (0.8770,0.3611) [circle,draw=black,fill=black,inner sep=0pt,minimum size=\rad] {};
\node (9) at (2.414,2.414) [circle,draw=black,fill=black,inner sep=0pt,minimum size=\rad] {};
\node (10) at (-2.414,2.414) [circle,draw=black,fill=white,inner sep=0pt,minimum size=\rad] {};
\node (11) at (-0.8770,0.3611) [circle,draw=black,fill=white,inner sep=0pt,minimum size=\rad] {};
\node (12) at (-0.9185,0.2795) [circle,draw=black,fill=black,inner sep=0pt,minimum size=\rad] {};
\node (13) at (-1.153,0.1647) [circle,draw=black,fill=black,inner sep=0pt,minimum size=\rad] {};
\node (14) at (-1.153,-0.1647) [circle,draw=black,fill=white,inner sep=0pt,minimum size=\rad] {};
\node (15) at (-0.4142,-0.4142) [circle,draw=black,fill=white,inner sep=0pt,minimum size=\rad] {};
\node (16) at (0.07107,-0.07107) [circle,draw=black,fill=white,inner sep=0pt,minimum size=\rad] {};
\node (17) at (0.07107,0.07107) [circle,draw=black,fill=black,inner sep=0pt,minimum size=\rad] {};
\node (18) at (-0.4142,0.4142) [circle,draw=black,fill=black,inner sep=0pt,minimum size=\rad] {};
\node (19) at (-0.9779,0.4027) [circle,draw=black,fill=black,inner sep=0pt,minimum size=\rad] {};
\node (20) at (-0.9982,0.3038) [circle,draw=black,fill=white,inner sep=0pt,minimum size=\rad] {};
\node (21) at (-0.8673,0.1239) [circle,draw=black,fill=white,inner sep=0pt,minimum size=\rad] {};
\node (22) at (-0.8673,-0.1239) [circle,draw=black,fill=black,inner sep=0pt,minimum size=\rad] {};
\node (23) at (-2.414,-2.414) [circle,draw=black,fill=black,inner sep=0pt,minimum size=\rad] {};
\node (24) at (2.414,-2.414) [circle,draw=black,fill=white,inner sep=0pt,minimum size=\rad] {};
\node (25) at (0.8673,-0.1239) [circle,draw=black,fill=white,inner sep=0pt,minimum size=\rad] {};
\node (26) at (0.8673,0.1239) [circle,draw=black,fill=black,inner sep=0pt,minimum size=\rad] {};
\node (27) at (0.9982,0.3038) [circle,draw=black,fill=black,inner sep=0pt,minimum size=\rad] {};
\node (28) at (0.9779,0.4027) [circle,draw=black,fill=white,inner sep=0pt,minimum size=\rad] {};

\draw[line width=1.2mm] (27)--(28)--(1)--(2)--(3);

\draw[line width=1.2mm] (6)--(7)--(8);

\draw[line width=1.2mm] (10)--(11)--(12);

\draw[line width=1.2mm] (14)--(15)--(16);

\draw[line width=1.2mm] (20)--(21)--(22);

\draw[line width=1.2mm] (24)--(25)--(26);

\draw[line width=2.5mm,white] (3)--(4)--(5)--(6);
\draw[line width=1.2mm] (3)--(4)--(5)--(6);

\draw[line width=2.5mm,white] (8)--(9)--(10);
\draw[line width=1.2mm] (8)--(9)--(10);

\draw[line width=2.5mm,white] (12)--(13)--(14);
\draw[line width=1.2mm]  (12)--(13)--(14);

\draw[line width=2.5mm,white] (16)--(17)--(18)--(19)--(20);
\draw[line width=1.2mm]  (16)--(17)--(18)--(19)--(20);

\draw[line width=2.5mm,white] (22)--(23)--(24);
\draw[line width=1.2mm]  (22)--(23)--(24);

\draw[line width=2.5mm,white] (26)--(27);
\draw[line width=1.2mm]  (26)--(27);
\end{tikzpicture} 
		\caption{An orthocubic representation of the knot $8_{19}$ with $28$ spheres (left) and its cubic diagram (right).}\label{fig:p323}
	\end{figure}

	In the nonalternating case, it is possible to construct orthocubic algebraic tangles with a total number of spheres strictly less than $4$ times the crossing number. The first non-trivial example that we have found satisfying this property is the Pretzel knot $P(3,-2,3)$, which corresponds to the knot $8_{19}$ in the Alexander-Briggs-Rolfsen notation. This knot is not alternating \cite{cromwell2004knots} and it admits an orthocubic necklace representation with $28$ spheres ($=\frac72\mathrm{cr}(8_{19})$, see Figure \ref{fig:p323}). However, it becomes more tricky to establish a relation with the crossing number in the non-alternating case since, in general, the crossing number is not the sum of the number of crossings of its rational factors.
	\medskip
	
\section{A new visualization of the slope of rational tangles}\label{sec:points}
The slope $p/q$ of a rational tangle $t_{p/q}$ has several geometric interpretations. For instance, it can be identified with the slope of the meridian of a solid torus that is the branched double covering of a rational tangle \cite{cromwell2004knots}. In this last section, we  present another geometric interpretation of the correspondance between rational tangles and rational numbers. We do so by relating the slope of a tangle with the slope of the line passing through the origin and the frist tangency point in the orthocubic Conway's construction. Astonishingly, this approach turns out to be helpful to find infinitely many primitive solutions of the Diophantine equation $x^4+y^4+z^4=2t^2$. 
\medskip

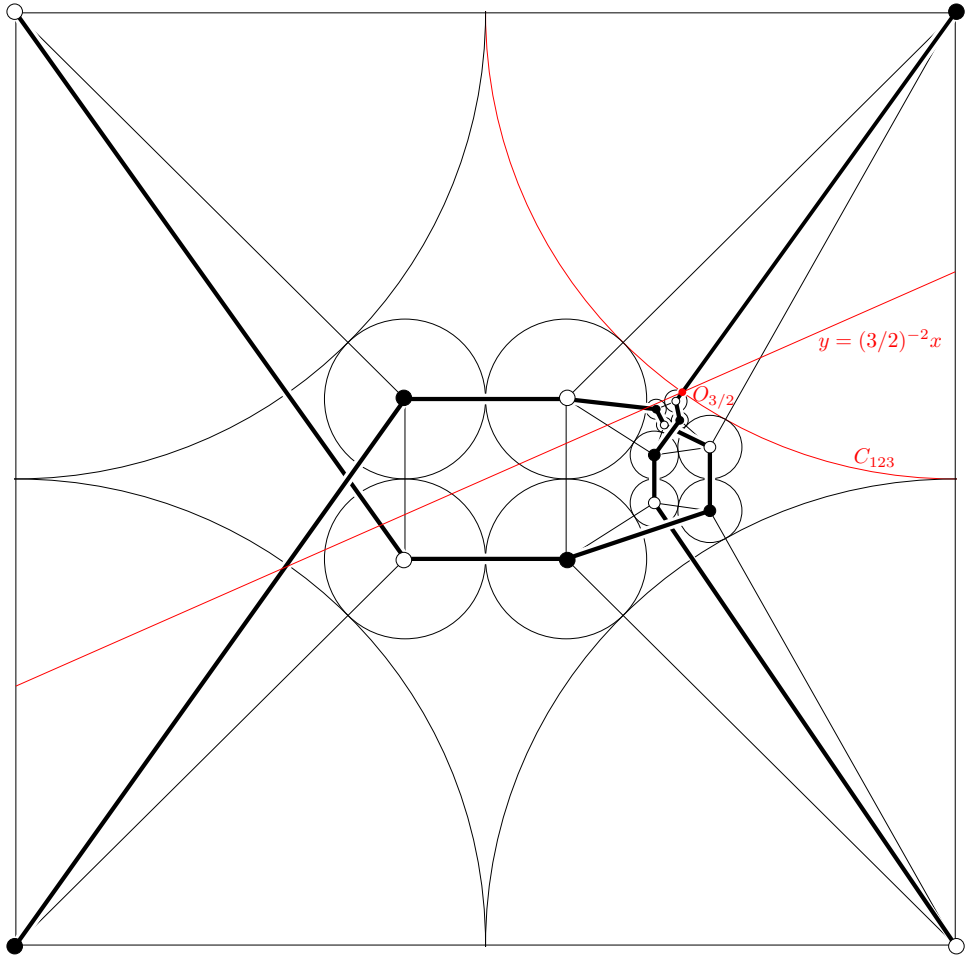
\begin{figure}[H]
					\centering
	\begin{tikzpicture}[scale=3.2]
		

		\begin{scope}

			\clip (-2.42,-2.42) rectangle (2.42,2.42);

			\draw[red] (2.414,2.414) circle (2.414cm) node at (2,.11) {$C_{123}$};
			\draw[very thin] (-2.414,2.414) circle (2.414cm);
			\draw[very thin] (2.414,-2.414) circle (2.414cm);
			\draw[very thin] (-2.414,-2.414) circle (2.414cm);
			\draw[very thin] (.414,.414) circle (.414cm);
			\draw[very thin] (-.414,.414) circle (.414cm);
			\draw[very thin] (.414,-.414) circle (.414cm);
			\draw[very thin] (-.414,-.414) circle (.414cm);
			
			\draw[very thin] (1.153,0.1647) circle (0.1647cm);
			\draw[very thin] (0.8673,0.1239) circle (0.1239cm);
			\draw[very thin] (1.153,-0.1647) circle (0.1647cm);
			\draw[very thin] (0.8673,-0.1239) circle (0.1239cm);
			\draw[very thin] (0.9779,0.4027) circle (0.05752cm);
			\draw[very thin] (0.8770,0.3611) circle (0.05159cm);
			\draw[very thin] (0.9982,0.3038) circle (0.04340cm);
			\draw[very thin] (0.9185,0.2795) circle (0.03993cm);

			\draw[thin] (2.414,2.414) -- (0.4142,0.4142);
			\draw[thin] (2.414,2.414) -- (2.414,-2.414);
			\draw[thin] (2.414,2.414) -- (-2.414,2.414);
			\draw[thin] (2.414,2.414) -- (1.15301,0.164716);
			\draw[thin] (2.414,2.414) -- (0.977867,0.402651);
			\draw[thin] (0.4142,0.4142) -- (0.4142,-0.4142);
			\draw[thin] (0.4142,0.4142) -- (-0.4142,0.4142);
			\draw[thin] (0.4142,0.4142) -- (0.867295,0.123899);
			\draw[thin] (0.4142,0.4142) -- (0.876977,0.361108);
			\draw[thin] (2.414,-2.414) -- (0.4142,-0.4142);
			\draw[thin] (2.414,-2.414) -- (-2.414,-2.414);
			\draw[thin] (2.414,-2.414) -- (1.15301,-0.164716);
			\draw[thin] (0.4142,-0.4142) -- (-0.4142,-0.4142);
			\draw[thin] (0.4142,-0.4142) -- (0.867295,-0.123899);
			\draw[thin] (-2.414,2.414) -- (-0.4142,0.4142);
			\draw[thin] (-2.414,2.414) -- (-2.414,-2.414);
			\draw[thin] (-0.4142,0.4142) -- (-0.4142,-0.4142);
			\draw[thin] (-2.414,-2.414) -- (-0.4142,-0.4142);
			\draw[thin] (1.15301,0.164716) -- (0.867295,0.123899);
			\draw[thin] (1.15301,0.164716) -- (1.15301,-0.164716);
			\draw[thin] (1.15301,0.164716) -- (0.998193,0.303798);
			\draw[thin] (0.867295,0.123899) -- (0.867295,-0.123899);
			\draw[thin] (0.867295,0.123899) -- (0.918471,0.279535);
			\draw[thin] (1.15301,-0.164716) -- (0.867295,-0.123899);
			\draw[thin] (0.977867,0.402651) -- (0.876977,0.361108);
			\draw[thin] (0.977867,0.402651) -- (0.998193,0.303798);
			\draw[thin] (0.876977,0.361108) -- (0.918471,0.279535);
			\draw[thin] (0.998193,0.303798) -- (0.918471,0.279535);

			\draw[line width=3pt,white] (2.414,-2.414)--(0.8673,-0.1239);
			\draw[line width=2pt,black] (2.414,-2.414)--(0.8673,-0.1239);

			\draw[line width=5pt,white] (-.414,.414)--(.414,.414)--(0.8770,0.3611)--(0.9185,0.2795)--(1.153,0.1647)--(1.153,-0.1647)--(.414,-.414)--(-.414,-.414)--(-2.414,2.414);
			\draw[line width=2pt] (-.414,.414)--(.414,.414)--(0.8770,0.3611)--(0.9185,0.2795)--(1.153,0.1647)--(1.153,-0.1647)--(.414,-.414)--(-.414,-.414)--(-2.414,2.414);
			
			\draw[line width=5pt,white] (-2.414,-2.414)--(-.414,.414);
			\draw[line width=2pt,black] (-2.414,-2.414)--(-.414,.414);
			
			\draw[line width=5pt,white] (0.9779,0.4027)--(2.414,2.414);
			\draw[line width=1.8pt] (0.9779,0.4027)--(2.414,2.414); 
			
			\draw[line width=5pt,white] (0.8673,-0.1239)--(0.8673,0.1239)--(0.9982,0.3038)--(0.9779,0.4027) ;
			\draw[line width=2pt,black](0.8673,-0.1239)--(0.8673,0.1239)--(0.9982,0.3038)--(0.9779,0.4027) ;

			\draw[red] (-2.41421,-1.07298) -- (2.41421,1.07298) node [pos=.92,below, inner sep=.4cm] {$y=(3/2)^{-2}x$};
					\fill[red] (1.01129, 0.449464)  circle (.02cm) node[xshift=.52cm,yshift=-.09cm] {$O_{3/2}$} ;
		\end{scope}
	\draw[fill=black]  (2.42,2.42) circle (.4mm);
	\draw[fill=white]  (2.42,-2.42) circle (.4mm);
	\draw[fill=white] (-2.42,2.42) circle (.4mm);
	\draw[fill=black] (-2.42,-2.42) circle (.4mm);
	
	\draw[fill=black]  (-0.42,0.42) circle (.4mm);
	\draw[fill=white]  (-0.42,-0.42) circle (.4mm);
	\draw[fill=white] (0.42,0.42) circle (.4mm);
	\draw[fill=black] (0.42,-0.42) circle (.4mm);
	
	\draw[fill=white] (1.153,0.1647) circle (.3mm);
	\draw[fill=black]  (0.8673,0.1239) circle (0.3mm);
	\draw[fill=black] (1.153,-0.1647) circle (.3mm);
	\draw[fill=white] (0.8673,-0.1239) circle (.3mm);
	\draw[fill=white] (0.9779,0.4027) circle (.2mm);
	\draw[fill=black] (0.8770,0.3611) circle (.2mm);
	\draw[fill=black] (0.9982,0.3038) circle (.2mm);
	\draw[fill=white] (0.9185,0.2795) circle (.2mm);
	\end{tikzpicture}  
						\vspace{-.3cm}
					\caption{The orthocubic point (red) of the rational tangle $t_{3/2}$ corresponding to the primitive solution $3^4+2^4+1^4=2\times7^2$.}\label{fig:fundamentals}
				\end{figure}

Let $p/q$ be a positive fraction with positive continued fraction expansion $[a_1,\cdots,a_n]$. We define the \textit{orthocubic point} $O_{p/q}$ of the rational tangle $t_{p/q}$ as the tangency point of the two circles in the cubic diagram of $t_\O(a_1,\cdots,a_n)$ corresponding to the \textit{first edge} of the orthocubic tangle. By \textit{first edge}, we mean the edge connecting the circle $C_{123}\in\C_{\{4,3\}}$  in the upper-right corner (see Figure \ref{fig:fundamentals}). We point out that this circle corresponds to the sphere $S_{4}\in\S_{\{3,3,4\}}$, which remains fixed under the orthocubic Conway's algorithm. On the other hand, the first edge is the projection of the image of the edge $(\overline 2,4)$, which is a nondiagonal edge. Since the orthocubic Conway's algorithm preserves the $z$-coloring, the first edge will be always a nondiagonal edge, so the orthocubic point is well-defined for every positive fraction. We can naturally extend the notion of orthocubic point to tangles with negative fractions  by applying a reflection through the plane $\{x=0\}$ to the whole setting. 

\begin{thm}\label{thm:slope}
For every $p/q\in \mathbb{Q}^{\pm}\cup\{\infty\}$, $O_{p/q}$ is the closest intersection to the origin of the line $y=\pm(p/q)^{-2}x$ with the circle $C_{\pm123}\in\C_{\{4,3\}}$.
\end{thm}
\begin{proof}
	
	It is enough to prove the positive case. Let $p\ge0 $ and $q\ge1$ be two coprime integers. We claim that 
	\begin{align}\label{eq:etacoord}
		\mathbf i(O_{p/q})=\begin{pmatrix}
			p^2\\q^2\\(p-q)^2\\\sqrt{2}(p^2-pq+q^2)
		\end{pmatrix}.
	\end{align}
	This would imply that $O_{p/q}=(\sqrt{2}(p^2-pq+q^2)-(p-q)^2)^{-1}(p^2,q^2)$. One can check that this is exactly the closest intersection to the origin of the line $\{y=(p/q)^{-2}x\}$ with the circle centred at $(1+\sqrt{2},1+\sqrt{2})$ and radius $(1+\sqrt{2})$, which is the boundary of the circle $C_{123}\in\C_{\{4,3\}}$.
	\medskip

	Let us prove the equality \eqref{eq:etacoord}.	The positiveness of $p$ and $q$ implies that we can find a positive continued fraction expansion $[a_1,\cdots,a_n]=p/q$  with $a_1\ge0$ and $a_i\ge1$ for every $1<i\le n$. Let $\boxed{t_{p/q}}$ the orthocubic tangle $t_\O[a_1,\ldots,a_n]$. Let $O_{p/q}$ and $O_\infty$ be the orthocubic points of $t_{p/q}$ and $t_\infty$, respectively. Now, by combining the definitions of the orthocubic operations $H_\O$ and $F_\O$, with the group morphism $\Gamma_{\{4,3\}}\longrightarrow\Gamma_{\{3,3,4\}} $ induced by the equivariant bijection $\phi_{\{3,3,4\}}^{\{4,3\}}$, and the definition of orthocubic rational tangles given in \eqref{eq!orthoConway}, we have that
	\begin{align*}
		\boxed{t_{p/q}}=H_\O^{a_1} F_\O\cdots H_\O^{a_n} F_\O\boxed{t_\infty}\Rightarrow
		O_{p/q}&=\mu_1^{a_1}r_{12}\cdots \mu_x^{a_n}r_{12}(O_\infty)\\
		&=(s_1r_{13})^{a_1}r_{12}\cdots (s_1r_{13})^{a_n}r_{12}(O_\infty)
	\end{align*}
	where $s_1$, $r_{13}$ and $r_{12}$ are the symmetries of $\Gamma_{\{4,3\}}$ described in Section \ref{subsec:orthoshifts}.  The matrices $\mathbf{S}_1$, $\mathbf{R}_{12}$, and $\mathbf{R}_{13}$ in $\mathrm{SL}_4(\mathbb R)$ representing these symmetries, and the inversive coordinates of $O_\infty$ can be computed by using the equations \eqref{eq:invcoordpoint} (with $\lambda=\sqrt{2}-1$) and \eqref{eq:invmatrix}, giving		
		\begin{align*}
			&&\mathbf{S}_1={\small\begin{pmatrix}
					-3 & 0 & 0 & 2 \sqrt{2} \\
					0 & 1 & 0 & 0 \\
					0 & 0 & 1 & 0 \\
					-2 \sqrt{2} & 0 & 0 & 3 \\
			\end{pmatrix}},&&
			\mathbf{R}_{13}={\small\begin{pmatrix}
					0 & 0 & 1 & 0 \\
					0 & 1 & 0 & 0 \\
					1 & 0 & 0 & 0 \\
					0 & 0 & 0 & 1 \\
			\end{pmatrix}},&& 
\mathbf{R}_{12}={\small\begin{pmatrix}
					0 & 1 & 0 & 0 \\
					1 & 0 & 0 & 0 \\
					0 & 0 & 1 & 0 \\
					0 & 0 & 0 & 1 \\
			\end{pmatrix}},					&&	\mathbf{i}(O_\infty)={\small\begin{pmatrix}
			1\\
			0\\
			1\\
			\sqrt2
		\end{pmatrix}}.
		\end{align*}
		
		Let $\mathbf{M}(k):=(\mathbf{S}_1\mathbf{R}_{13})^k\mathbf{R}_{12}$. By induction on $k$, we have
		\begin{align*}
			\mathbf{M}(k)=\begin{pmatrix}
				0 & 1-k^2 & -k (k+2) & \sqrt{2} k (k+1) \\
				1 & 0 & 0 & 0 \\
				0 & -k(k-2) & 1-k^2 & \sqrt{2} k (k-1) \\
				0 & -\sqrt{2} k(k-1) & -\sqrt{2} k (k+1) & 2 k^2+1 \\
			\end{pmatrix}
		\end{align*}
		We prove the equality \eqref{eq:etacoord} by induction on the number of coefficients $n$ in the fraction expansion of $p/q$. For $n=1$ (that is $p=a_1$ and $q=1$) we have
		\begin{align*}
			\mathbf i(O_{a_1})=\mathbf{M}(a_1){\small\begin{pmatrix}
					1\\
					0\\
					1\\
					\sqrt{2}
			\end{pmatrix}}={\small\begin{pmatrix}
					a_1^2\\
					1\\
					(a_1-1)^2\\
					\sqrt{2}(a_1^2-a_1+1)
			\end{pmatrix}}.
		\end{align*}
		\noindent We suppose that equality \eqref{eq:etacoord} holds true for $n-1\ge1$. Let $r/s=a_2+\frac{1}{\cdots+\frac{1}{a_n}}$. Then,
		\begin{align*}
			\mathbf i(O_{p/q})=\mathbf{M}(a_1)\cdots\mathbf{M}(a_n)\begin{pmatrix}
				1\\
				0\\
				1\\
				\sqrt{2}
			\end{pmatrix}&=\mathbf{M}(a_1)\begin{pmatrix}
				r^2\\
				s^2\\
				(r-s)^2\\
				\sqrt{2}(r^2-rs+s)
			\end{pmatrix}
		&=\begin{pmatrix}
				(ra_1+s)^2\\
				r^2\\
				(ra_1+s-r)^2\\
				\sqrt{2}((ra_1+s)^2-r(ra_1+s)+r^2)
			\end{pmatrix}
		\end{align*}
		We finally notice that $\frac{ra_1+s}{r}=a_1+\frac{1}{r/s}=a_1+\frac{1}{a_2+\frac{1}{\cdots+\frac{1}{a_n}}}=p/q$ and therefore, equality \eqref{eq:etacoord} holds.
	\end{proof}

				\begin{cor}\label{cor:diaphortho}The Diophantine equation $x^4+y^4+z^4=2t^2$
					has an infinite number of primitive solutions.
				\end{cor}
				\begin{proof}
				By combining equations \eqref{eq:invprod}, \eqref{eq:invprodpoint} with the inversive coordinates of the orthocubic point of every rational tangle $t_{p/q}$ given in \eqref{eq:etacoord}, we obtain the following parametrization generating infinite primitive solutions 
					\begin{align}
						x=p,&&y=q,&&z=p-q,&&t=p^2-pq+q^2.
					\end{align}
					\vspace{-.3cm}
				\end{proof}

	The parametrization of primitive solutions provided above can be found through other classic methods. However, we find interesting its geometric interpretation as tangency points of the orthocubic tangle contruction (see Figure \ref{fig:fundamentals}). This connection arises from in the fact that the full symmetry group $\Gamma_{\{4,3\}}$, which generates the cubic crystallographic packing, is an integral arithmetic subgroup of the orthogonal group for the quadratic form $x^2+y^2+z^2-2t^2$ \cite{stange2015bianchi,KontorovichNakamura}. Points of tangency in this packing represent integral isotropic points for the quadratic form. According to Theorem \ref{thm:slope}, the family of orthocubic points corresponding to positive fractions provides an infinite family of points of tangency where $x$, $y$, $z$ are squares.	We hope and expect that this approach will be useful to find solutions to other type of Diophantine equations.
\vspace{-.5cm}

	\printbibliography[
title={References}
] 

\end{document}